\documentclass[11pt,reqno]{amsart}
\usepackage{fixltx2e}                    
\usepackage[osf]{mathpazo}  
\usepackage[utf8]{inputenc}            
\usepackage{amsmath}                     
\usepackage{amssymb, latexsym, stmaryrd, amsthm, dsfont, amsfonts, amsbsy,amsthm, amsmath, mathrsfs}            
\usepackage{mathtools}                   
\usepackage{bm}                          
\usepackage{enumerate}                   
\usepackage{verbatim}                    
\usepackage{url}   
\usepackage{lscape}                      

\usepackage{microtype}  
\usepackage[all, knot]{xy}
        \xyoption{arc} 
        \xyoption{web}                 
\makeatletter                            
\def\MT@register@subst@font{\MT@exp@one@n\MT@in@clist\font@name\MT@font@list
 \ifMT@inlist@\else\xdef\MT@font@list{\MT@font@list\font@name,}\fi}
\makeatother

\usepackage[pdftex,bookmarks,bookmarksnumbered,linktocpage,   
         colorlinks,linkcolor=blue,citecolor=blue]{hyperref}

\newcommand{\bit}{\begin{itemize}}    
\newcommand{\eit}{\end{itemize}}
\newcommand{\ben}{\begin{enumerate}}
\newcommand{\een}{\end{enumerate}}
\newcommand{\benormal}{\ben[\normalfont 1.]}   

\let\enormal\een
\newcommand{\benroman}{\ben[\normalfont (i)]}  
\let\eroman\een

\newcommand{\bde}{\begin{description}}
\newcommand{\ede}{\end{description}}

\newcommand{\?}{\ensuremath{\mkern0.4\thinmuskip}}   

\let\leq=\leqslant
\let\geq=\geqslant
\let\Box=\square                            

\let\epsilon=\varepsilon
\let\Lambda\varLambda
\let\Gamma\varGamma
\let\Delta\varDelta
\let\Lambda\varLambda
\let\Omega\varOmega
\let\Theta\varTheta
\let\Xi\varXi
\let\Pi\varPi
\let\Sigma\varSigma

\let\class=\mathsf                              
\let\oper=\mathbb                               

\bmdefine{\A}{A}                                
\bmdefine{\2}{2}
\bmdefine{\B}{B}
\bmdefine{\D}{D}
\bmdefine{\M}{M}                                
\bmdefine{\LLL}{L}                              
\bmdefine{\Fm}{Fm}                              
\bmdefine{\zerou}{[0{,}1]}  
\bmdefine{\T}{T}                                

\newcommand{\VVV}{\oper{V}}                     
\newcommand{\QQQ}{\oper{Q}}

\newcommand{\HHH}{\oper{H}}

\newcommand{\PPP}{\oper{P}}
\newcommand{\PSD}{\oper{P}_{\!\textsc{sd}}^{}}
\newcommand{\PPU}{\oper{P}_{\!\textsc{u}}^{}}

\newcommand{\SSS}{\oper{S}}
\newcommand{\III}{\oper{I}}

\newcommand{\Con}{\mathrm{Con}\?}                            
      
\bmdefine{\boldstar}{\mathchoice{\textstyle*}{\textstyle*}{\textstyle*}{\scriptstyle*}}

\bmdefine{\btau}{\tau}                                  
\bmdefine{\brho}{\rho}                                  

\newcommand{\semeq}{\mathrel{=\joinrel\mathrel\vert\mkern0mu\mathrel\vert\joinrel=}}  

\bmdefine{\leibniz}{\Omega}        

\bmdefine{\frege}{\Lambda}         

\makeatletter
\newcommand{\tarskidsp}{\mathord%
   {\m@th\raisebox{0pt}[0pt][0pt]{$\stackrel%
   {\raisebox{-2.7pt}[0ex][0pt]{$\displaystyle \,\?\thicksim$}}%
   {\displaystyle\leibniz}$}}}
\newcommand{\tarskitxt}{\mathord%
   {\m@th\raisebox{0pt}[0pt][0pt]{$\stackrel%
   {\raisebox{-2.7pt}[0ex][0pt]{$\,\?\thicksim$}}{\displaystyle\leibniz}$}}}
\newcommand{\tarskiscr}{\mathord%
   {{\m@th\raisebox{0pt}[0pt][0pt]{$\stackrel%
   {\raisebox{-2.4pt}[0ex][0pt]{$\scriptstyle \,\?\thicksim$}}%
   {\scriptstyle\leibniz}$}}}}
\newcommand{\tarskiscrscr}{\mathord%
   {{\m@th\raisebox{0pt}[0pt][0pt]{$\stackrel%
   {\raisebox{-2pt}[0ex][0pt]{$\scriptscriptstyle \,\?\thicksim$}}%
   {\scriptscriptstyle\leibniz}$}}}}
\newcommand{\tarski}{\@ifnextchar ^ %
   {\mathchoice{\tarskidsp\kern-.07em}{\tarskitxt\kern-.07em}%
   {\tarskiscr\kern-.07em}{\tarskiscrscr\kern-.07em}}%
   {\mathchoice{\tarskidsp}{\tarskitxt}{\tarskiscr}{\tarskiscrscr}}}
\makeatother

\theoremstyle{theorem}
\newtheorem{Theorem}{Theorem}[section]
\newtheorem{Lemma}[Theorem]{Lemma}
\newtheorem{Corollary}[Theorem]{Corollary}
\newtheorem{Fact}[Theorem]{Fact}

\theoremstyle{definition}
\newtheorem{law}[Theorem]{Definition}
\newtheorem{exa}[Theorem]{Example}

\theoremstyle{remark}

\newtheorem{problem}{\bf Problem}

\newcommand{\C}{\boldsymbol{C}} 

\begin{document}
\title[Varieties of positive modal algebras and structural completeness]{Varieties of positive modal algebras and structural completeness}

\author{Tommaso Moraschini}
\email{moraschini@cs.cas.cz}
\address{Institute of Computer Science of the Czech Academy of Science, Prague, Czech Republic}
\date{\today}
\renewcommand{\thefootnote}{\fnsymbol{footnote}} 
\renewcommand{\thefootnote}{\arabic{footnote}} 
\maketitle


\begin{abstract}
Positive modal algebras are the $\langle \land, \lor, \Diamond, \Box, 0, 1\rangle$-subreducts of modal algebras.\ We prove that the variety of positive S4-algebras is not locally finite.\ On the other hand, the free one-generated positive S4-algebra is shown to be finite.\ Moreover, we describe the bottom part of the lattice of varieties of positive S4-algebras. Building on this, we characterize (passively, hereditarily) structurally complete varieties of positive K4-algebras.
\end{abstract}

\section{Introduction}

In recent years, subclassical fragments of modal logics have attracted some attention because of their balance between expressive power and lower computational complexity \cite{Beklemishev14a,Bekle18,Nguyen2000b,Nguyen05a}.\ Among these subclassical fragments, a special role is played by the \textit{positive modal logic} PML \cite{Du95}, i.e.\ the $\langle \land, \lor, \Box, \Diamond, 0, 1 \rangle$-fragment of the local consequence of the normal modal logic ${\bf K}$ \cite{BlRiVe01,ChZa97,HuCr96,Kr99}. The interest of PML and its extensions comes from the fact that, despite their restricted signature, they preserve some mathematically desiderable features typical of normal modal logics, including Sahlqvist theory and a well-behaved Priestley-style duality \cite{CeJa97,CeJa99a,GeNaVe05}. For similar reasons, PML and their generalizations were studied both from the perspective of algebraic logic \cite{Ja02} and from that of coalgebra \cite{BaKuVe15,Pa03}.

When suitably reformulated as a Gentzen system $\mathcal{PML}$ as in \cite{CeJa97}, the logic PML becomes \textit{algebraizable} in the sense of \cite{Ra06,ReV95-p} (see the Appendix for a more detailed discussion). In particular, the algebraic counterpart of $\mathcal{PML}$ is the class of the so-called \textit{positive modal algebras} \cite{Du95}, i.e.\ positive subreducts of classical modal algebras, which turn out to form a variety (or, equivalently, by Birkhoff’s theorem, an equational class). As a consequence of algebraizability, there is a dual lattice isomorphism between the lattice of axiomatic extensions of $\mathcal{PML}$ and that of subvarieties of positive modal algebras. Thus axiomatic extensions of $\mathcal{PML}$ can be studied through the lenses of varieties of positive modal algebras, which in turn are amenable to the methods of universal algebra \cite{Be11g,BuSa00,McMcTa87}. Driven by this observation, in this paper we initiate the study of the lattice of subvarieties of positive modal algebras. Since its structure is largely unknown, we focus on the more definite goal of describing its \textit{structurally complete} members. 

For a general background on structural completeness we refer the reader to \cite{RaxxNJ,Ry97}. Even if structural completeness is usually understood as a property of logics and not, strictly speaking, of classes of algebras, we find it convenient to introduce this and some related concepts in purely algebraic terms as in \cite{Be91b}. Let $\class{K}$ be a quasi-variety, and consider a quasi-equation
\[
\Phi \coloneqq \varphi_{1} \thickapprox \psi_{1}\& \dots \& \varphi_{n} \thickapprox \psi_{n} \to \varphi \thickapprox \psi.
\]
The quasi-equation $\Phi$ is \textit{active} in $\class{K}$ when there is a substitution $\sigma$ such that $\class{K}\vDash \sigma \varphi_{i} \thickapprox \sigma\psi_{i}$ for every $i \leq n$. The quasi-equation $\Phi$ is \textit{passive} in $\class{K}$ if it is not active in $\class{K}$. The quasi-equation $\Phi$ is \textit{admissible} in $\class{K}$ if for all substitutions $\sigma$:
\[
\text{if }\class{K} \vDash \sigma \varphi_{i} \thickapprox \sigma \psi_{i}\text{ for every }i \leq n\text{, then }\class{K}\vDash \sigma \varphi \thickapprox \sigma \psi.
\]
The quasi-equation $\Phi$ is \textit{valid} in $\class{K}$ if $\class{K}\vDash \Phi$. Accordingly, we say that
\benormal
\item $\class{K}$ is \textit{actively structurally complete} (ASC) if every active admissible quasi-equation is valid.
\item $\class{K}$ is \textit{passively structurally complete} (PSC) if every passive quasi-equation is valid.\footnote{Observe that passive quasi-equations are vacuously admissible.}
\item $\class{K}$ is \textit{structurally complete} (SC) if every admissible quasi-equation is valid.
\item $\class{K}$ is \textit{hereditarily structurally complete} (HSC) if every subquasi-variety of $\class{K}$ is SC.
\enormal

If $\vdash$ is a finitary algebraizable logic or Gentzen system \cite{BP89,Ra06,ReV95-p} with equivalent algebraic semantics a quasi-variety $\mathsf{K}$, then $\mathsf{K}$ is structurally complete exactly when the finite admissible rules of $\vdash$ are derivable in the usual sense. This is the case, for instance, of the Gentzen system $\mathcal{PML}$, which is finitary and algebraized by the variety of positive modal algebras.

The interest of variants of structural completeness is partially due to their influence on the structure of the lattices of subquasi-varieties. For instance, if a quasi-variety $\mathsf{K}$ is HSC, then its lattice of subquasi-varieties happens to be distributive \cite{Gorbunov1976,Go98a}. Moreover, if $\mathsf{K}$ is PSC, then all its subquasi-varieties have the joint embedding property \cite{MorRafWan18} and, therefore, are generated as quasi-varieties by a single algebra \cite{Ma71}. On the other hand, variants of structural completeness are interesting also from a purely logical perspective. For instance, if a SC quasi-variety $\mathsf{K}$ has the finite model property (i.e.\ $\mathsf{K}$ and the class of its finite members generate the same \textit{variety}), then $\mathsf{K}$ has the strong finite model property as well (i.e.\ it is generated \textit{as a quasi-variety} by its finite members), as observed in \cite{OlRaAl08}. Finally, ASC is known to be related to projective and exact unification \cite{DzSt201x}.

 In \cite{Be91b,DzSt201x,MeRo13a,Wro09} the following characterizations of the various kinds of structural completeness are obtained.

\begin{Theorem}\label{Thm:Bergman} Let $\class{K}$ be a variety.
\
\benormal
\item $\class{K}$ is ASC if and only if $\A \times \Fm_{\class{K}}(\omega) \in \QQQ (\Fm_{\class{K}}(\omega))$ for every subdirectly irreducible algebra $\A \in \class{K}$. If there is a constant symbol in the language, then we can replace ``$\A \times \Fm_{\class{K}}(\omega) \in \QQQ (\Fm_{\class{K}}(\omega))$'' by ``$\A \times \Fm_{\class{K}}(0) \in \QQQ (\Fm_{\class{K}}(\omega))$'' in this statement.
\item $\class{K}$ is a PSC variety if and only if every positive existential sentence is either true in all non-trivial members of $\class{K}$ or false in all non-trivial members of $\class{K}$. 
\item $\class{K}$ is SC if and only if $\A \in \SSS\PPU (\Fm_{\class{K}}(\omega))$ for every subdirectly irreducible algebra $\A \in \class{K}$.
\item $\class{K}$ is HSC if and only if every subquasi-variety of $\class{K}$ is a variety.
\enormal
\end{Theorem}

Structural completeness and its variants have been the subject of intense research in the setting of modal and Heyting algebras \cite{Ry97}, yielding an array of deep results.\ To name only a few of them, HSC varieties of K4-algebras, and of Heyting algebras have been fully described in \cite{Ryb95,Citkin78a}. Moreover, explicit bases for the admissible rules of intuitionistic logic (equiv.\ of quasi-equations admissible in the variety of Heyting algebras), and of some prominent modal systems have been provided in \cite{Ie01,Jerabek05c}. Unfortunately, the picture of structural completeness tends to change dramatically when the set of connectives is altered (see for instance \cite{RaftSwiry16}). Accordingly, we should not expect, at least in principle, that results on admissibility and structural completeness valid in modal algebras could be extended directly to positive modal algebras.

Our results on structural completeness are confined to classes of \textit{positive K4-algebras}, i.e.\ positive subreducts of K4-algebras.\ In this context, we show that there are exactly three non-trivial structurally complete varieties of positive K4-algebras (Theorem \ref{Thm:StructuralCompleteness++}). Moreover, these are also the unique hereditarily structurally complete varieties of positive K4-algebras. This contrasts with the full-signature case, since there are infinitely (countably) many hereditarily structurally complete varieties of K4-algebras \cite{Ryb95}.\ Moreover, we characterize passively structurally complete varieties of positive K4-algebras (Theorem \ref{Thm:PSC}) and prove that there are infinitely many of them (Example \ref{Ex:PSC}). 

To establish these results, we rely on two observations concerning \textit{positive S4-algebras}, i.e.\ positive subreducts of S4-algebras, which may be of some interest on their own. On the one hand, we show that, even if the variety of positive S4-algebras is not locally finite (Corollary \ref{Cor:Non-locfinite}), its free one-generated algebra is finite (Theorem \ref{Thm:Free-one-genalgebra}). This contrasts with the full-signature case, since already the one-generated free S4-algebra is infinite \cite{Blo77}. On the other hand, we rely on a description of the bottom part of the lattice of subvarieties of positive S4-algebras up to height $\leq 4$ (Theorem \ref{Thm:Height4}), represented pictorially in Figure \ref{Fig:Varieties}. From a logical point of view, this result provides a partial description of  the top part of the lattice of axiomatic extensions of the Gentzen system $\mathcal{PML}$. 

We hope that this work may stimulate the further investigation of varieties of positive modal algebras and of admissibility in positive modal logic, whose theory is at the moment underdeveloped. In particular, it would be interesting to find an explicit base for the admissible quasi-equations of the variety of positive K4-algebras, complementing for instance the results of \cite{Jerabek05c}.

\section{Preliminaries}

For a general background on universal algebra we refer the reader to \cite{Be11g,BuSa00,McMcTa87}. We denote by $\III, \HHH, \SSS, \PPP, \PSD, \PPU$ the class operators of closure under isomorphism, homomorphic images, subalgebras, direct products, subdirect products and ultraproducts respectively. Given an algebra $\A$,  we denote by $\Con\A$ its congruence lattice and by $\textup{Id}_{\A}$ the identity relation on $A$. Then $\A$ is \textit{subdirectly irreducible} when $\textup{Id}_{\A}$ is completely meet-irreducible in $\Con\A$, i.e.\ when the poset $\Con\A \smallsetminus \{ \textup{Id}_{\A} \}$ has a least element $\theta$. In this case, the congruence $\theta$ is called the \textit{monolith} of $\A$. The monolith of $\A$ is always a principal congruence, in fact $\theta = \textup{Cg}(a, b)$ for every $\langle a, b\rangle \in \theta$ such that $a \ne b$. Similarly, an algebra $\A$ is \textit{finitely subdirectly irreducible} if $\textup{Id}_{\A}$ is meet-irreducible in $\Con\A$. Given a class of algebras $\class{K}$, we denote by $\class{K}_{si}$ the class of its subdirectly irreducible members. An algebra $\A$ is \textit{simple} when $\Con \A$ has exactly two-elements.   

A \textit{variety} is a class of algebras axiomatized by equations or, equivalently, a class of algebras closed under $\HHH,\SSS$ and $\PPP$. A \textit{quasi-variety} is a class of algebras axiomatized by quasi-equations or, equivalently, a class of algebras closed under $\III,\SSS, \PPP$ and $\PPU$. Given a class of algebras $\class{K}$, we denote by $\VVV(\class{K})$ and $\QQQ(\class{K})$ respectively the least variety and quasi-variety containing $\class{K}$. It is well known that $\VVV(\class{K}) = \HHH\SSS\PPP(\class{K})$ and $\QQQ(\class{K})= \III\SSS\PPP\PPU(\class{K})$. An  algebra $\A$ in a quasi-variety $\class{K}$ is \textit{subdirectly irreducible relative} to $\class{K}$ if the poset $\langle \{ \theta \in \Con\A : \A / \theta \in \class{K} \text{ and }\theta \ne \textup{Id}_{\A} \}, \subseteq \rangle$ has a minimum element. Every algebra in a quasi-variety $\class{K}$ is a subdirect product of algebras subdirectly irreducible relative to $\class{K}$. If $\class{K}$ is a class of algebras, then the algebras subdirectly irreducible relative to $\QQQ(\class{K})$ belong to $\III\SSS\PPU(\class{K})$ \cite{CzDz}. Given a quasi-variety $\class{K}$ and $n \in \omega$, we denote by $\Fm_{\class{K}}(x_{1}, \dots, x_{n})$ the free $n$-generated algebra over $\class{K}$. It is well known that the elements of $\Fm_{\class{K}}(x_{1}, \dots, x_{n})$ are congruence classes of terms, which we denote by $\llbracket \varphi\rrbracket, \llbracket \psi \rrbracket$ etc. The free countably-generated algebra over $\class{K}$ is denoted by $\Fm_{\class{K}}(\omega)$.

Let $\class{K}$ be a class of algebras. An algebra $\A \in \class{K}$ is \textit{projective} in $\class{K}$ if for every $\B, \C \in \class{K}$ and homomorphisms $f \colon \B \to \C$ and $g \colon \A \to \C$ where $f$ is surjective, there is a homomorphism $h \colon \A \to \B$ such that $g= f\circ h$. An algebra $\A$ is a \textit{retract} of an algebra $\B$ if there are homomorphisms $f \colon \A \to \B$ and $g \colon \B \to \A$ such that the composition $g \circ f$ is the identity map on $A$. It is well known that the projective algebras in a variety are exactly the retracts of free algebras.

Consider an algebra $\A$ and a sublanguage $\mathscr{L}$ of the language of $\A$. The $\mathscr{L}$-\textit{reduct} of $\A$ is the algebra $\langle A, \{f^{\A} : f \in \mathscr{L}\} \rangle$. An algebra $\B$ in language $\mathscr{L}$ is a $\mathscr{L}$-\textit{subreduct} of $\A$ if it can be embedded into the $\mathscr{L}$-reduct of $\A$. If $\class{K}$ is a quasi-variety and $\mathscr{L}$ is a sublanguage of the language of $\class{K}$, then the class of $\mathscr{L}$-subreducts of algebras in $\class{K}$ coincides with the quasi-variety generated by the $\mathscr{L}$-reducts of algebras in $\class{K}$ \cite{Ma71}.

A variety $\class{K}$ has the \textit{congruence extension property} (CEP) if for every $\A, \B \in \class{K}$ such that $\B \in \SSS(\A)$, if $\theta$ is a congruence of $\B$, then there is a congruence $\phi$ of $\A$ such that $\phi \cap B^{2} = \theta$. A variety $\class{K}$ has \textit{equationally definable principal congruences} (EDPC) \cite{KoPi80} if there is a finite set of equation $\Phi(x, y, z, v)$ such that for every $\A \in \class{K}$ and $a, b, c, d \in A$ we have:
\[
\langle a, b\rangle \in \textup{Cg}(c, d) \Longleftrightarrow \A \vDash \Phi(a, b, c, d).
\]
It is easy to see that EDPC implies the CEP, however the converse implication does not hold in general.\ It is well known that both the variety of (bounded) distributive lattices and the variety of K4-algebras have EDPC. More precisely, given a (bounded) distributive lattice $\A$ and $a, b, c, d \in A$, we have that
\begin{equation}\label{Eq:EDPCDL}
\langle c, d\rangle \in \textup{Cg}(a, b) \Longleftrightarrow ( c \land a \land b = d \land a \land b \text{ and }c \lor a \lor b = d \lor a \lor b).
\end{equation}
Moreover, given a K4-algebra $\A$ and $a, b, c, d \in A$, we have that
\begin{equation}\label{Eq:EDPCK4}
\langle a, b\rangle \in \textup{Cg}(c, d) \Longleftrightarrow \A \vDash (c \leftrightarrow d) \land \Box (c \leftrightarrow d) \leq a \leftrightarrow b.
\end{equation}

Finally, we will rely on the following easy observation, which is essentially  \cite[Theorem 7.7]{RaxxNJ}.

\begin{Lemma}\label{Cor:PSC_simple}
If $\class{K}$ is a PSC variety with a constant symbol, its free $0$-generated algebra is either simple or trivial.
\end{Lemma}

\section{Algebras and frames}

The algebraic study of the positive fragment of the normal modal logic ${\bf K}$ was begun by Dunn in \cite{Du95}, where the following algebraic models were introduced:

\begin{law}
A \textit{positive modal algebra} is a structure $\A = \langle A, \land, \lor, \Box, \Diamond, 0, 1\rangle$ where $\langle A, \land, \lor, 0, 1\rangle$ is a bounded distributive lattice such that $\Box 1 =1$ and $\Diamond 0 = 0$ and
\begin{align*}
\Box (a \land b) &= \Box a \land \Box b\\
\Diamond ( a \lor b) &= \Diamond a \lor \Diamond b\\
\Box a \land \Diamond b &\leq \Diamond (a \land b)\\
\Box (a \lor b) & \leq \Box a \lor \Diamond b
\end{align*}
for every $a, b \in A$. We denote by $\class{PMA}$ the variety of positive modal algebras. The class $\class{PMA}$ can be endowed with the structure of a category, whose arrows are homomorphisms, and which we will denote also by $\class{PMA}$ since no confusion shall occur.
\end{law}

It is worth remarking that the above equations were considered earlier by Johnstone in the study of a localic version of the Vietoris construction \cite{Jo82,PTJo85}.


A \textit{Priestley space} is a triple $\langle X, \leq, \tau\rangle$, where $\langle X, \leq \rangle$ is a poset and $\langle X, \tau\rangle$ is a compact topological space, which satisfies the following additional condition: for every $x, y \in X$ such that $x \nleq y$ there is a clopen upset $U$ such that $x \in U$ and $y \notin U$. Priestley duality states that the category of Priestley spaces, endowed with continuous order-preserving maps, is dually equivalent to the category of bounded distributive lattices with homomorphisms \cite{Pr70,Pr72,DaPr02}. Building on Priestley duality, Celani and Jansana presented a topological duality for positive modal algebras \cite{CeJa99a}. We will briefly review it, since it will be needed later on. Given a set $X$ and a binary relation $R$ on it, we let $\Box_{R}, \Diamond_{R} \colon \mathcal{P}(X) \to \mathcal{P}(X)$ be the maps defined for every $V \subseteq X$ as follows:
\begin{align*}
\Box_{R}V &\coloneqq \{ x \in X : \text{ if }\langle x, y\rangle \in R\text{, then }y \in V \}\\
\Diamond_{R} V &\coloneqq \{ x \in X: \text{ there is }y\in X \text{ such that }\langle x, y\rangle \in R\text{ and }y \in V \}.
\end{align*}

\begin{law}\label{Def:K+space}
A $\class{K}^{+}$\textit{-space} is a structure $\langle X, \leq, R, \tau\rangle$ where $\langle X, \leq, \tau\rangle$ is a Priestley space and $R$ is a binary relation on $X$ such that:
\benormal
\item $R = ( R \circ \leq ) \cap (R \circ \leq^{-1})$.
\item The clopen upsets of $\tau$ are closed under the operations $\Box_{R}$ and $\Diamond_{R}$.
\item The set $\{ y \in X : \langle x, y \rangle \in R \}$ is topologically closed for every $x \in X$.
\enormal
\end{law}

\begin{law}
Let $\boldsymbol{X}=\langle X, \leq^{\boldsymbol{X}}, R^{\boldsymbol{X}}, \tau^{\boldsymbol{X}}\rangle$ and $\boldsymbol{Y}=\langle Y, \leq^{\boldsymbol{Y}}, R^{\boldsymbol{Y}}, \tau^{\boldsymbol{Y}}\rangle$ be $\class{K}^{+}$-spaces. A \textit{p-morphism} from $\boldsymbol{X}$ to $\boldsymbol{Y}$ is a continuous map $f \colon X \to Y$, that preserves $\leq$ and $R$, and satisfies the following condition for every $x, y \in X$:
\begin{align*}
\text{if }\langle f(x), y \rangle \in R^{\boldsymbol{Y}}&\text{, then there are }z, v \in X\text{ such that }\\
&\?\?\?\?\?\langle x, z\rangle, \langle x, v\rangle \in R^{\boldsymbol{X}}\text{ and }f(z)\leq y \leq f(v).
\end{align*}
We denote by $\class{K}^{+}$ the category of $\class{K}^{+}$ spaces endowed with p-morphisms.
\end{law}

Positive modal algebras and $\class{K}^{+}$-spaces are related as follows. Consider $\A \in \class{PMA}$ and let $\textup{Pr}(\A)$ be the collection of prime filters of $\A$. Then we define a relation $R_{\A} \subseteq \textup{Pr}(\A) \times \textup{Pr}(\A)$ as follows:
\[
\langle F, G \rangle \in R_{\A} \Longleftrightarrow \Box^{-1}(F) \subseteq G \subseteq \Diamond^{-1}(F)
\]
for every $F, G\in \textup{Pr}(\A)$. Moreover, for every $a \in A$ we set
\[
\varphi(a) \coloneqq \{ F \in \textup{Pr}(\A) : a \in F \}.
\]
Finally, we let $\tau$ be the topology on $\textup{Pr}(\A)$ generated by the following subbasis:
\[
\{ \varphi(a) : a \in A \} \cup \{ \varphi(a)^{c} : a \in A \}.
\]
It turns out that the structure
\[
\A^{\ast} \coloneqq \langle \textup{Pr}(\A), \subseteq, R_{\A}, \tau\rangle
\]
is a $\class{K}^{+}$-space. Moreover, if $f \colon \A \to \B$ is a homomorphism between positive modal algebras, then the map $f^{\ast} \colon \textup{Pr}(\B) \to \textup{Pr}(\A)$, defined for every $F \in \textup{Pr}(\B)$ as $f^{\ast}(F) \coloneqq f^{-1}[F]$, is a p-morphism from $\B^{\ast}$ to $\A^{\ast}$. It turns out that the application $^{\ast} \colon \class{PMA} \to \class{K}^{+}$ is a contravariant functor.

Conversely, consider a $\class{K}^{+}$-space $\boldsymbol{X}$ and let $\textup{Up}(\boldsymbol{X})$ be the collection of clopen upsets of $\langle X, \leq, \tau\rangle$. It turns out that the structure
\[
\boldsymbol{X}_{\ast} \coloneqq \langle \textup{Up}(X), \cap, \cup, \Box_{R}, \Diamond_{R}, \emptyset, X\rangle
\]
is a positive modal algebra. Moreover, if $f \colon \boldsymbol{X} \to \boldsymbol{Y}$ is a p-morphism between $\class{K}^{+}$-spaces, then the map $f_{\ast} \colon \textup{Up}(Y) \to \textup{Up}(X)$, defined for every $V \in \textup{Up}(Y)$ as $f_{\ast}(V) \coloneqq f^{-1}[V]$, is a homomorphism from $\boldsymbol{Y}_{\ast}$ to $\boldsymbol{X}_{\ast}$. The application $_{\ast} \colon \class{K}^{+} \to \class{PMA}$ is a contravariant functor as well. The relation between positive modal algebras and $\class{K}^{+}$-spaces can be formulated as follows \cite[pag. 700]{CeJa99a}:

\begin{Theorem}[Celani and Jansana]\label{Thm:CategoryEquivalence}
The functors $^{\ast} \colon \class{PMA} \longleftrightarrow \class{K}^{+} \colon _{\ast}$ from a dual category equivalence.
\end{Theorem}

The above result implies that every positive modal algebra $\A$ is isomorphic to the algebra of clopen upsets of its dual $\class{K}^{+}$-space. More precisely, we have the following representation theorem, which is essentially \cite[Theorem 8.1]{Du95}, but see also \cite[Theorem 2.2]{CeJa99a}.
\begin{Theorem}[Celani, Dunn and Jansana]\label{Thm:Representation}
Consider $\A \in \class{PMA}$. The map $\kappa \colon \A \to (\A^{\ast})_{\ast}$ defined by the rule
\[
a \longmapsto \{ F \in \textup{Pr}(\A) : a \in F \}
\]
is an isomorphism.
\end{Theorem}

A useful consequence of this representation is the following correspondence result:

\begin{Corollary}\label{Cor:Correspondence}
Consider $\A \in \class{PMA}$.
\benormal
\item $R_{\A}$ is reflexive if and only if $\Box a \leq a \leq \Diamond a$ for every $a \in A$.
\item $R_{\A}$ is transitive if and only if $\Box a \leq \Box \Box a$ and $\Diamond \Diamond a \leq \Diamond a$ for every $a \in A$.
\enormal
\end{Corollary}

We will focus on the following varieties of positive modal algebras:

\begin{law}
Let $\A \in \class{PMA}$.
\benormal
\item $\A$ is a \textit{positive K4-algebra} if for every $a \in A$:
\[
\Box a \leq \Box \Box a \text{ and }\Diamond \Diamond a \leq \Diamond a.
\]
\item $\A$ is a \textit{positive S4-algebra} if for every $a \in A$:
\[
\Box \Box a = \Box a \leq a \leq \Diamond a = \Diamond \Diamond a.
\]
\enormal
We denote respectively by $\class{PK4}$ and $\class{PS4}$ the varieties of positive K4 and S4-algebras. It is clear that $\class{PS4} \subseteq \class{PK4}$.
\end{law}

From Theorem  \ref{Thm:CategoryEquivalence} and Corollary \ref{Cor:Correspondence} it follows that $\class{PK4}$ (resp.\ $\class{PS4}$) is dually equivalent to the full subcategory $\class{K}^{+}$, whose objects are $\class{K}^{+}$-spaces whose relation $R$ is transitive (resp.\ transitive and reflexive). This dual equivalence is witnessed by the restriction of the functors described above. The following result motivates the name of positive K4 and S4-algebras:

\begin{Theorem}\label{Thm:Subreducts}
$\class{PMA}$, $\class{PKA}$ and $\class{PS4}$ are respectively the classes of $\langle \land, \lor, \Box, \Diamond, 0, 1\rangle$ subreducts of modal, K4 and S4-algebras.
\end{Theorem}
\begin{proof}
Consider $\A \in \class{PMA}$ and define an algebra
\[
\mathcal{M}(\A) \coloneqq \langle \mathcal{P}(\textup{Pr}(\A)), \cap, \cup, \lnot, \Box_{R_{\A}}, 0, 1\rangle
\]
where $\lnot$ is the set-theoretic complement. It is very easy to see that $\mathcal{M}(\A)$ is a modal algebra, and that  that map $\kappa \colon \A \to \mathcal{M}(\A)$, defined in Theorem \ref{Thm:Representation}, is an embedding that respects $\langle \land, \lor, \Box, \Diamond, 0, 1\rangle$. Thus $\A$ is a positive subreduct of a modal algebra.

Now consider $\A \in \class{PK4}$. By Corollary \ref{Cor:Correspondence} we know that $R_{\A}$ is transitive. In particular, this implies that $\mathcal{M}(\A)$ is a K4-algebra. Thus the map $\kappa \colon \A \to \mathcal{M}(\A)$ embeds $\A$ into a K4-algebra preserving $\langle \land, \lor, \Box, \Diamond, 0, 1\rangle$. We conclude that $\A$ is a positive subreduct of a K4-algebra. On the other hand, the equations defining $\class{PK4}$ hold in K4-algebras. Thus we conclude that $\class{PK4}$ is the class of positive subreducts of K4-algebras. The case of $\class{PS4}$ is handled similarly.
\end{proof}

\section{Well-connected algebras}

Let $\A$ be a modal algebra. A lattice filter $F$ of $\A$ is \textit{open} if $a \in F$ implies $\Box a \in F$ for every $a \in A$. The set of open filters of $\A$, when ordered under the inclusion relation, forms a lattice $\textup{Op}(\A)$. It is well known \cite{Kr99} that the lattices $\textup{Op}(\A)$ and $\Con\A$ are isomorphic under the map defined by the following rule:
\[
F \longmapsto \theta_{F} \coloneqq \{ \langle a, b \rangle \in A^{2} : a \to b, b \to a \in F \}.
\]
An analogous correspondence holds between the congruences of $\A$ and its ideals closed under $\Diamond$. The correspondence between open filters and congruences implies that a modal algebra $\A$ is subdirectly irreducible when the poset $\langle \textup{Op}(\A) \smallsetminus \{ \{ 1 \} \}, \subseteq \rangle$ has a minimum element \cite{Rauten1979bk}. In some cases this condition can be shifted from subsets of the universe of $\A$ (such as open filters) to elements of $\A$. This happens for example in S4-algebras.\footnote{Similar results hold for K4-algebras as well.} An S4-algebra $\A$ is \textit{well-connected} if $\Box a \lor \Box b = 1$ implies $a=1$ or $b=1$ for every $a, b \in A$ \cite{McKT44}. The correspondence between open filters and congruences implies that an S4-algebra is well-connected if and only if it is finitely subdirectly irreducible. Moreover, an S4-algebra $\A$ is simple if and only if $\Box a = 0$ for every $a \ne 1$ (see \cite{Kr99} if necessary). The notion of well-connection can be adapted to positive S4-algebras as follows:
\begin{law}
 A positive S4-algebra $\A$ is \textit{well-connected} when for every $a, b \in A$:
\benormal
\item If $\Box a \lor \Box b = 1$, then $a=1$ or $b=1$.
\item If $\Diamond a \land \Diamond b = 0$, then $a= 0$ or $b = 0$.
\enormal
\end{law}
It is worth observing that free positive S4-algebras are well-connected. This follows immediately from the well-known fact that free S4-algebras are well-connected or, equivalently, from the fact that the modal system $\boldsymbol{S4}$ has the modal disjunction property.

The results in this section show that certain aspects of finitely subdirectly irreducible and simple algebras are preserved from S4-algebras to their positive subreducts.
Given a positive modal algebra $\A$, let $\mathcal{M}(\A)$ be the modal algebra defined the proof of Theorem \ref{Thm:Subreducts}. 

\begin{Theorem}\label{Thm:FSI-transfers}
 If $\A \in \class{PS4}$ is finitely subdirectly irreducible, then so is $\mathcal{M}(\A)$.
\end{Theorem}

\begin{proof}
Consider a finitely subdirectly irreducible algebra $\A \in \class{PS4}$. The fact that $\A$ is a positive subreduct of an S4-algebra, together with the correspondence between open filters and congruences typical of modal algebras, easily implies that $\A$ is well-connected.

Since $\A$ is well-connected, the following sets are, respectively, a proper filter and a proper ideal of $\A$:
\begin{align*}
F &\coloneqq \{ a \in A : \Diamond b \leq a \text{ for some }b \in A\smallsetminus \{0\}\}\\
I &\coloneqq \{ a \in A : a \leq \Box b \text{ for some }b \in A\smallsetminus \{1\}\}.
\end{align*}

We claim that $F \cap I = \emptyset$. Suppose the contrary in view of a contradiction. Then there is $a \in F \cap I$, i.e.\ there are $b, c \in A$ such that
\[
0 < \Diamond b \leq a \leq \Box c < 1.
\]
Consider the embedding $\kappa \colon \A \to \mathcal{M}(\A)$ defined in the proof of Theorem \ref{Thm:Subreducts}. Let $F$ and $G$ be respectively the upsets of $\lnot \kappa(\Diamond b)$ and $\kappa(\Box c)$ in $\mathcal{M}(\A)$. From the fact that $\Diamond b$ and $\Box c$ are respectively fixed points of $\Diamond$ and $\Box$, it follows that $F$ and $G$ are open filters. Since $\Diamond b \leq \Box c$, we have that $F \cap G = \{ 1 \}$. This means that $\theta_{F} \cap \theta_{G} = \textup{Id}_{\mathcal{M}(\A)}$. Define
\begin{align*}
\phi_{F} &\coloneqq \{ \langle x, y \rangle \in A^{2} : \langle \kappa(x), \kappa(y)\rangle \in \theta_{F}\}\\
\phi_{G} &\coloneqq \{ \langle x, y \rangle \in A^{2} : \langle \kappa(x), \kappa(y)\rangle \in \theta_{G}\}.
\end{align*}
Clearly we have that $\phi_{F} \cap \phi_{G} = \textup{Id}_{\A}$. Moreover, neither $\phi_{F}$ not $\phi_{G}$ are the identity relation, since $\langle \Diamond b, 0\rangle \in \phi_{F}$ and $\langle \Box c, 1 \rangle \in \phi_{G}$. But this contradicts the fact that $\A$ is finitely subdirectly irreducible, establishing the claim.

From the claim it follows that $F \cap I = \emptyset$. From the Prime Filter Theorem it follows that there is a prime filter $F^{\ast}$ extending $F$ and such that $F^{\ast} \cap I = \emptyset$. In particular, we have that for every $a \in A$,
\begin{align*}
(\Box a \in F^{\ast} \Longrightarrow a = 1) \text{ and }(\Diamond a \notin F^{\ast} \Longrightarrow a = 0).
\end{align*}
This implies that for every $G \in \textup{Pr}(\A)$,
\begin{equation}\label{Eq:FSI-transfer1}
\Box^{-1}F^{\ast} \subseteq G \subseteq \Diamond^{-1}F^{\ast}, \text{ for every }G \in \textup{Pr}(\A).
\end{equation}

Suppose in view of a contradiction that $\mathcal{M}(\A)$ is not subdirectly irreducible, i.e.\ that it is not well-connected. Then there are two elements $X, Y \in \mathcal{M}(\A) \smallsetminus \{1\}$ such that $\Box X \lor \Box Y = 1$. Observe that
\[
X, Y \in \mathcal{P}(\textup{Pr}(\A))\smallsetminus \{ \textup{Pr}(\A) \} \text{ and }\Box X \cup \Box Y = \textup{Pr}(\A).
\]
In particular, we have that $F^{\ast} \in \Box X \cup \Box Y$. We can assume w.l.o.g.\ that $F^{\ast} \in \Box X$. From the definition of the operation $\Box$ in $\mathcal{M}(\A)$ it follows that for every $G \in \textup{Pr}(\A)$,
\[
\Box^{-1}F^{\ast} \subseteq G \subseteq \Diamond^{-1}F^{\ast} \Longrightarrow G \in X.
\]
Together with (\ref{Eq:FSI-transfer1}), this implies that $X = \textup{Pr}(\A)$ which is a contradiction. Hence we conclude that $\mathcal{M}(\A)$ is finitely subdirectly irreducible.
\end{proof}

\begin{Corollary}\label{Cor:WellConnected}
Consider $\A \in \class{PS4}$ finitely subdirectly irreducible.  $\A$ is well-connected and for every $a \in A \smallsetminus \{ 0, 1 \}$ either $\Box a < a$ or $a < \Diamond a$.
\end{Corollary}

However the converse of Theorem \ref{Thm:FSI-transfers} does not hold in general, as shown in the next example. The failure of the converse is due to the fact that the correspondence between open filters and congruences is lost, when we move from S4-algebras to their positive subreducts.

\begin{exa}
It is natural to wonder whether properties such as \textit{being well-connected} or \textit{being finitely subdirectly irreducible} are preserved or reflected by the passage from a positive S4-algebra $\A$ to the S4-algebra $\mathcal{M}(\A)$. It turns out that the situation is as follows:
\benroman
\item If $\A$ is finitely subdirectly irreducible, then so is $\mathcal{M}(\A)$.
\item If $\mathcal{M}(\A)$ is well-connected, then so is $\A$.
\item Even if $\A$ is well-connected, $\mathcal{M}(\A)$ may fail to be so.
\item Even if $\mathcal{M}(\A)$ is finitely subdirectly irreducible, $\A$ may fail to be so.
\eroman

Observe that (i) amounts to Theorem \ref{Thm:FSI-transfers}, while (ii) is clear.

(iii): Observe that in a positive S4-algebra we always have that $0$ and $1$ are fixed points of the modal operations. Consider the positive S4-algebra $\A$, whose lattice reduct is the three element chain $0 < a  < 1$ and such that $\Box a = \Diamond a = 1$. Clearly $\A$ is well-connected. However, $\mathcal{M}(\A)$ is the four element modal algebra where the modal operations are interpreted as identity maps. In particular, $\mathcal{M}(\A)$ is not well-connected.

(iv): Consider the positive S4-algebra $\A$, whose lattice reduct is the five element chain $0 < a < b < c < 1$ and such that $\Box e = 0$ and $\Diamond e = 1$ for every $e \in \{ a, b, c \}$. Consider the following equivalence relations, described through their blocks:
\[
\theta \coloneqq \{ \{ 0 \}, \{a, b \}, \{ c \}, \{ 1 \} \} \text{ and }\phi \coloneqq \{ \{ 0 \}, \{a \}, \{b, c \}, \{ 1 \} \}.
\]
Clearly $\theta, \phi \in \Con\A$ and $\theta \cap \phi = \textup{Id}_{\A}$. Thus the identity relation is not meet-irreducible in $\Con \A$ and, therefore, $\A$ is not finitely subdirectly irreducible. However, it is not hard to see that $\mathcal{M}(\A)$ is a simple S4-algebra.
\qed
\end{exa}

The failure of the correspondence between open filters and congruences makes it hard to imagine a concrete description of subdirectly irreducible positive S4-algebras.\footnote{It is worth remarking that subdirectly irreducible positive modal algebras have been characterized by means of topological properties of their corresponding $\class{K}^{+}$-spaces \cite{Ce06a}. However, it is not straightforward to apply this characterization to the study of problems such as the classifications of varieties of positive modal algebras.}\ 	Additional difficulties in the description of subdirectly irreducible positive modal algebras arise from the fact that the CEP (and, therefore, also EDPC) fails in $\class{PS4}$ as shown in Example \ref{Exa:CEP}. Nevertheless it is easy to obtain a complete description of all simple positive K4-algebras. To this end, let $\boldsymbol{B}_{2}$ be the unique two-element positive modal algebra such that $\Diamond 1 = 0$ and $\Box 0 = 1$. 

\begin{Lemma}\label{Lem:Simple}
Let $\A$ be a non-trivial member of $\class{PK4}$. $\A$ is simple if and only if either $\A = \B_{2}$ or the following conditions hold:
\benroman
\item For every $a \in A$:
\[
\Box a\coloneqq \left\{ \begin{array}{ll}
 1 & \text{if $a = 1$}\\
 0 & \text{otherwise}
  \end{array} \right.    
\quad  \Diamond a \coloneqq \left\{ \begin{array}{ll}
 0 & \text{if $a=0$}\\
 1 & \text{otherwise.}
  \end{array} \right.  
\]
\item For every $a, b \in A$ such that $0 < a < b < 1$, there is $c \in A \smallsetminus \{ 0, 1 \}$ such that
\[
\text{either }(a \leq c \text{ and }b \lor c = 1) \text{ or }(c \leq b \text{ and }a \land c = 0).
\]
\eroman
\end{Lemma}

\begin{proof}
We begin by proving the ``only if'' part. Suppose that $\A$ is simple and not isomorphic to $\B_{2}$. We claim that $\Box 0 \ne 1$. To prove this suppose towards a contradiction that $\Box 0 = 1$. Then clearly by monotonicity $\Box a = 1$ for all $a \in A$. From the fact that $\Box 0 =1$, it follows that $\Diamond 1 = 0$. Together with the fact that $\A$ is not isomorphic to $\B_{2}$, we obtain that there is an element $a \in A$ such that $0 < a < 1$. Then consider the embedding $\kappa \colon \A \to \mathcal{M}(\A)$ defined in the proof of Theorem \ref{Thm:Subreducts}. Let $F$ be the upset of $\kappa(a)$ in $\mathcal{M}(\A)$. From the fact that $\Box a = 1$, it follows that $F$ is an open filter. We have that
\[
\langle \kappa(a), 1 \rangle \in \theta_{F} \text{ and }\langle \kappa(a), 0\rangle \notin \theta_{F}.
\]
We define $\phi \coloneqq \{ \langle b, c \rangle \in A^{2} : \langle \kappa(b), \kappa(c)\rangle \in \theta_{F}\}$.  Clearly $\phi$ is a congruence of $\A$, that is neither the identity nor the total relation. Thus $\kappa(\A)$ is not simple, which is a contradiction. This establishes our claim. 

Now suppose towards a contradiction that there is an element $a \in A \smallsetminus \{ 1 \}$ such that $\Box a \ne 0$. Then consider the embedding $\kappa \colon \A \to \mathcal{M}(\A)$. Let $F$ be the principal filter generated by $\kappa(\Box a)$ in $\mathcal{M}(\A)$. We have that
\[
\langle \kappa(\Box a), 1\rangle \in \theta_{F} \text{ and }\langle \kappa(\Box a), 0 \rangle \notin \theta_{F}.
\]
Again we define $\phi \coloneqq \{ \langle b, c \rangle \in A^{2} : \langle \kappa(b), \kappa(c)\rangle \in \theta_{F}\}$. Clearly $\phi$ is a congruence of $\A$, that is not the total relation on $A$, since $\langle \Box a, 0 \rangle \notin \phi$. Since $\A$ is simple, this implies that $\phi$ is the identity relation on $A$. In  particular, this means that $\Box a =1$. From our claim it follows that $0 \ne a$. Then consider the upset $F$ of $\kappa(a)$ in $\mathcal{M}(\A)$. From the fact that $\Box a = 1$, it follows that $F$ is an open filter. We have that
\[
\langle \kappa(a), 0\rangle \notin \theta_{F}\text{ and }\langle \kappa(a), 1 \rangle \in \theta_{F}.
\]
We define $\phi \coloneqq \{ \langle b, c \rangle \in A^{2} : \langle \kappa(b), \kappa(c)\rangle \in \theta_{F}\}$. Clearly $\phi$ is a congruence of $\A$, that is neither the identity (since $a \ne 1$) not the total relation. Thus $\A$ is not simple, which is a contradiction. We conclude $\Box a = 0$ for all $a \ne 1$.  A similar argument shows that $\Diamond a = 1$ for every $a \ne 0$. The only difference from the one described above is that this time one need to rely on the correspondence between the congruences of $\mathcal{M}(\A)$ and its ideals closed under $\Diamond$.

Now, suppose towards a contradiction that there are $a, b \in A$ such that $0 < a < b < 1$ and for which condition (ii) of the statement fails. Consider the congruence $\theta$ of the bounded lattice $\langle A, \land, \lor, 0, 1 \rangle$ generated by the pair $\langle a, b \rangle$. Clearly $\theta$ is not the identity relation. Applying EDPC for distributive lattices as in (\ref{Eq:EDPCDL}), it is easy to see that $\{ 0 \}$ and $\{ 1 \}$ are blocks of $\theta$. Thus $\theta$ is not the total relation. Now we prove that $\theta$ is compatible w.r.t.\ $\Box$ and $\Diamond$. Consider a pair $\langle c, d \rangle \in \theta$ such that $c \ne d$. From the fact that $\{ 0 \}$ and $\{ 1 \}$ are blocks of $\theta$, it follows that $0 < c, d < 1$. As we showed, this means that $\Box c = \Box d = 0$ and $\Diamond c = \Diamond d = 1$. We conclude that $\theta$ is a congruence of $\A$, which contradicts the fact that $\A$ is simple. This establishes that condition (ii) in the statement holds.

Then we turn to prove the ``if'' part. Observe that $\B_{2}$ is simple on cardinality grounds. Then suppose that $\A$ satisfies conditions (i) and (ii) in the statement. Suppose towards a contradiction that $\A$ is not simple. Then there are two different elements $a, b \in A$ such that $\textup{Cg}(a, b)$ is not the total relation. We can assume w.l.o.g.\ that $a < b$. Observe that $0 < a$. Suppose towards a contradiction that $a= 0$. Then $b < 1$, since $\textup{Cg}(a, b)$ is not the total relation. From condition (i) we know that $\Diamond a = 0$ and $\Diamond b = 1$. But this implies that $\textup{Cg}(a, b)$ is the total relation, which is false. This establishes that $0 < a$. A similar argument shows that $b < 1$. Hence we have that $0 < a < b < 1$ and we can apply condition (ii) in the statement. This condition, together with EDPC for distributive lattices as in (\ref{Eq:EDPCDL}), implies that there is $c \in A \smallsetminus \{ 0, 1 \}$ such that either the pair $\langle c, 1\rangle$ or the pair $\langle c, 0\rangle$ belongs to $\textup{Cg}(a, b)$. Together with the fact that $\Diamond c = 1$ and $\Box c = 0$, this implies that $\langle 0, 1 \rangle \in \textup{Cg}(a, b)$. But this means that $\textup{Cg}(a, b)$ is the total congruence, which is a contradiction.
\end{proof}

As we promised, the next example shows that $\class{PS4}$ does not have the CEP and, therefore, it does not have EDPC either. This contrasts with the well-known fact that the variety of S4-algebras has EDPC.

\begin{exa}\label{Exa:CEP}
Consider any Boolean lattice $\langle A, \land, \lor, 0, 1 \rangle$ of at least $8$ elements. Then equip it with modal operations $\Box$ and $\Diamond$ defined as in condition (i) of Lemma \ref{Lem:Simple}. Let $\A$ be the resulting positive S4-algebra. Consider any chain chain $0 < a < b < 1$ in $\langle A, \leq\rangle$. Observe that setting $c \coloneqq a \lor \lnot b$, where $\lnot b$ is the Boolean complement of $b$, we have that
\[
b \lor c = 1 \text{ and }a \leq c.
\]
Hence $\A$ satisfies condition (ii) of Lemma \ref{Lem:Simple}. We conclude that $\A$ is a simple algebra. Now consider any chain $0 < a < b < 1$ in $\langle \A, \leq \rangle$. Observe that it is the universe of a subalgebra $\B$ of $\A$. It is easy to see that $\B$ is not simple. Therefore $\A$ is a simple algebra with a non-simple and non-trivial subalgebra. We conclude that $\class{PS4}$ does not have the CEP.
\qed
\end{exa}

\section{Free one-generated positive S4-algebra}

Observe that the free one-generated positive K4-algebra is infinite. This can be easily seen by observing that it contains the infinite ascending chain
\[
\llbracket \Box x \rrbracket < \llbracket \Box^{2} x \rrbracket < \llbracket \Box^{3} x \rrbracket < \dots
\]
In this section we will prove that the situation changes for positive S4-algebras, whose free one-generated algebra turns out to be finite. This fact  contrasts with the full-signature case, since it is well known that the free one-generated S4-algebra is infinite, see for instance \cite{Blo77}.

For the sake of simplicity, let us introduce a way of describing pictorially positive S4-algebras. Let $\A$ be a positive S4-algebras. Then the structure of $\A$ is uniquely determined by its Hasse diagram and the fixed points of the modal operations. This is because for every $a \in A$, the element $\Box a$ (resp.\ $\Diamond a$) is the greatest fixed point of $\Box$ below $a$ (resp.\ smallest fixed point of $\Diamond$ above $a$). Thus a way of describing the structure of $\A \in \class{PS4}$ is to depict its Hasse diagram and to mark with a $\Box$ (resp.\ with a $\Diamond$) the fixed points of $\Box$ (resp.\ of $\Diamond$). We omit these marks for $0$ and $1$, which are always fixed points of the modal operations, since $\A$ is a positive modal algebras. An example of this way of representing positive S4-algebras is offered in Figure \ref{Fig:Free}. 
\begin{figure}
\[
\xymatrix@R=16pt @C=20pt @!0{
&&&& *-{1 \bullet\ \ \?  }\ar@{-}[d] &&\\
&&&& *-{\Diamond \bullet \ \ \?\?\? } \ar@{-}[d]&&\\
&&&& *-{\bullet}\ar@{-}[dr]\ar@{-}[dl] &&\\
&&&*-{\bullet}\ar@{-}[d] \ar@{-}[dl] \ar@{-}[dr] && *-{\  \ \bullet \Diamond} \ar@{-}[dl]&\\
&& *-{\bullet} \ar@{-}[d]  \ar@{-}[dr] & *-{\bullet} \ar@{-}[d]  \ar@{-}[dl]\ar@{-}[dr]  & *-{\bullet} \ar@{-}[d]  \ar@{-}[dl]\ar@{-}[dr] &&\\
&& *-{\bullet}\ar@{-}[dl]\ar@{-}[dr] & *-{\bullet} \ar@{-}[d]\ar@{-}[dr]& *-{\bullet}\ar@{-}[dl]\ar@{-}[dr] & *-{\bullet}\ar@{-}[d]\ar@{-}[dl]&\\
& *-{\bullet} \ar@{-}[dl]\ar@{-}[dr]&& *-{\bullet} \ar@{-}[dl]\ar@{-}[dr]& *-{\bullet} \ar@{-}[d]\ar@{-}[dr]& *-{\bullet} \ar@{-}[dl]\ar@{-}[dr]&\\
*-{a \bullet \?\?\?\?\?\? } \ar@{-}[dr]&& *-{\bullet} \ar@{-}[dr]\ar@{-}[dl]&& *-{\bullet} \ar@{-}[d]\ar@{-}[dr]\ar@{-}[dl]& *-{ \ \ \bullet \Diamond} \ar@{-}[dl]& *-{\ \ \bullet \Box}\ar@{-}[dl]\\
& *-{\bullet}\ar@{-}[dr] && *-{\bullet}\ar@{-}[d]\ar@{-}[dr]\ar@{-}[dl] & *-{\bullet} \ar@{-}[dl]\ar@{-}[dr]& *-{\bullet}\ar@{-}[d]\ar@{-}[dl] &\\
&& *-{\bullet} \ar@{-}[d]\ar@{-}[dr]& *-{\bullet}\ar@{-}[dl]\ar@{-}[dr] & *-{\bullet} \ar@{-}[d]\ar@{-}[dl]& *-{\bullet}\ar@{-}[dl] &\\
&& *-{\bullet}\ar@{-}[dr] & *-{\bullet} \ar@{-}[d]& *-{\bullet}\ar@{-}[dl]\ar@{-}[dr] &&\\
&&& *-{\bullet}\ar@{-}[dr] && *-{\ \ \bullet \Box} \ar@{-}[dl]&\\
&&&& *-{\bullet}\ar@{-}[dd]&&\\
&&&& *-{ \Box \bullet \ \ \?\?\?}\ar@{-}[d]&&\\
&&&& *-{0 \bullet \ \ \?}&&
}
\]
\caption{The free one-generated positive S4-algebra}\label{Fig:Free}
\end{figure}
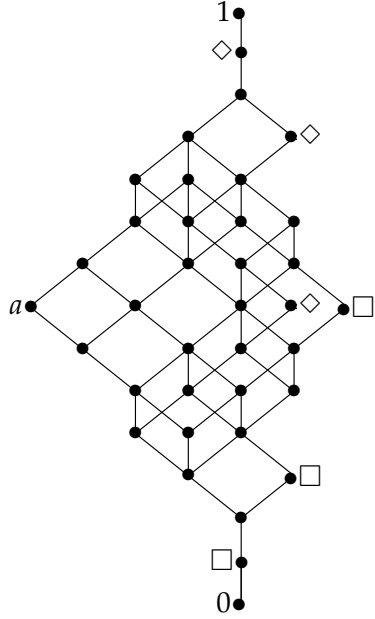
\begin{Theorem}\label{Thm:Free-one-genalgebra}
The algebra $\A$ depicted in Figure \ref{Fig:Free} is the free one-generated positive S4-algebra.
\end{Theorem}

\begin{proof}
Consider the terms
\begin{equation}
\Sigma \coloneqq \{ x, \Box x, \Diamond \Box x, \Box \Diamond \Box x, \Diamond x, \Box \Diamond x, \Box \Diamond \Box x \}.
\end{equation}
We will make use of the following easy observation:
\begin{Fact}\label{Fact:Free}
For every positive S4-algebra $\B$ and $b \in B$, the elements $\{ \varphi^{\B}(b) : \varphi \in \Sigma \} \subseteq B$ form the following subposet of $\langle B, \leq \rangle$, where the lines indicate only the order relation and must not be interpreted as referring to any description of meets and joins, and whose elements are not necessarily different one from the other:
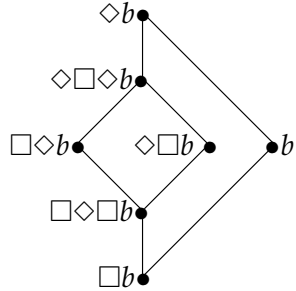
\[
\xymatrix@R=25pt @C=25pt @!0{
& *-{\Diamond b \bullet \ \ \ \ \?\?}\ar@{-}[d]\ar@{-}[ddrr] &&\\
& *-{\Diamond \Box \Diamond b \bullet \ \ \ \ \ \ \ \ \ \ \?\?\? }\ar@{-}[dr]\ar@{-}[dl] &&\\
*-{\Box \Diamond b \bullet  \ \ \ \ \ \?\?\?\?\?\?\?\?\?} \ar@{-}[dr]&& *-{\Diamond \Box b \  \bullet \ \ \ \ \ \ \?\?\?\?\?\?\?\?\? }\ar@{-}[dl] & *-{\?\?\?\?\? \bullet b}\ar@{-}[ddll]\\
& *-{\Box \Diamond \Box b \bullet  \ \ \ \ \ \ \ \ \ \ \?\?\? }\ar@{-}[d] &&\\
& *-{\Box b \bullet  \ \ \ \ \?\?} &&
}
\]
\end{Fact}

Now, let $\A$ be the algebra depicted in Figure \ref{Fig:Free}. It is an easy exercise to check that the algebra $\A$ is in fact a positive S4-algebra. Moreover, the bounded lattice reduct of $\A$ is obtained as follows. First we consider the free $7$-generated bounded distributive lattice $\C$, whose generators are the terms in $\Sigma$. Second we form the quotient of $\C$ under the congruence generated by following the set
\[
\Gamma \coloneqq \{ \langle \varphi, \psi \rangle \in \Sigma^{2} : \varphi(b) \leq \psi(b) \text{ in the above diagram}\}.
\]
The algebra $\C / \textup{Cg}(\Gamma)$ obtained in this way is exactly the bounded lattice reduct of $\A$. This can be checked mechanically, e.g.\ using the Universal Algebra Calculator \cite{UAC11}. In particular, for every element $c \in A$ there is a bounded lattice lattice term $t_{c}(y_{1}, \dots, y_{7})$ such that
\[
c = t_{c}^{\A}(a, \Box a, \Diamond \Box a, \Box \Diamond \Box a, \Diamond a, \Box \Diamond a, \Box \Diamond \Box a).
\]
We can assume w.l.o.g.\ that
\[
t_{0} = 0, t_{1} = 1\text{ and }t_{a} = y_{1}, t_{\Box a} = y_{2}, \dots, t_{\Box \Diamond \Box a} = y_{7}.
\]

In order to prove that $\A$ is the free one-generated positive S4-algebra with free generator $a$, it will be enough to show that $\A$ enjoys the following universal property: for every $\B \in \class{PS4}$ and $b \in B$, there is a unique homomorphism $f \colon \A \to \B$ such that $f(a) = b$. The uniqueness of the homomorphism $f$ is ensured by the fact that $\A$ is generated by $a$. Then it only remains to prove its existence. To this end, consider the map $f \colon A \to B$ defined by the following rule:
\[
f(c)\coloneqq t_{c}^{\B}(b, \Box b, \Diamond \Box b, \Box \Diamond \Box b, \Diamond b, \Box \Diamond b, \Box \Diamond \Box b).
\]
To prove that $f$ is a homomorphism of bounded lattices, we reason as follows. Clearly $f$ preserves $0$ and $1$, since $t_{0} = 0$ and $t_{1} = 0$. Moreover, the lattice reduct of $\A$ is freely-generated by $a, \Box a, \Diamond \Box a, \Box \Diamond \Box a, \Diamond a, \Box \Diamond a, \Box \Diamond \Box a$ w.r.t.\ to the inequalities in $\Gamma$, together with Fact \ref{Fact:Free}, implies that $f$ is a lattice homomorphism as well.

It only remains to prove that $f$ preserves the modal operations. We will detail only the fact that $f$ preserves $\Box$. Consider the following partition $\theta$ of $\A$:
\[
\{ 0 \}, [\Box a, a], [\Box \Diamond a, \Diamond a], [ \Box \Diamond \Box a, a  \lor \Diamond \Box a],  \{ 1 \}.
\]
Observe that if two elements $c, d \in A$ belong to the same block of $\theta$, then $\Box^{\A}c = \Box^{\A}d$. Since $f$ preserves $0$ and $1$, we have that
\[
f(\Box 0) = f(0) = 0 = \Box 0 \text{ and }f(\Box 1) = f(1) = 1 = \Box 1.
\]
Then, suppose that $c \in  [\Box a, a] \cup [\Box \Diamond a, \Diamond a] \cup [ \Box \Diamond \Box a, a  \lor \Diamond \Box a]$. We have that following cases:
\benormal
\item $c \in [ \Box \Diamond \Box a, a  \lor \Diamond \Box a]$.
\item $c \in  [\Box a, a]$.
\item $c \in [\Box \Diamond a, \Diamond a]$.
\enormal
1. Observe that
\[
\class{PS4} \vDash \Box ( x  \lor \Diamond \Box x ) \leq \Box x \lor \Diamond \Diamond \Box x = \Diamond \Box x.
\]
Together with the monotonicity of $\Box$, this implies that
\[
\Box ( x  \lor \Diamond \Box x ) = \Box\Box ( x  \lor \Diamond \Box x ) \leq  \Box \Diamond \Box x.
\]
Finally, from the fact that $\Diamond \Box x  \leq x  \lor \Diamond \Box x$, it follows that $\Box \Diamond \Box x \leq \Box ( x  \lor \Diamond \Box x )$. We conclude that
\begin{equation}\label{Eq:Free_trick1}
\class{PS4} \vDash \Box ( x  \lor \Diamond \Box x ) \thickapprox  \Box \Diamond \Box x.
\end{equation}
Since $c \in [ \Box \Diamond \Box a, a  \lor \Diamond \Box a]$ and $f$ is order preserving, we have that:
\[
\Box \Diamond \Box b = f(\Box \Diamond \Box a) \leq f(c) \leq f(a  \lor \Diamond \Box a) = b \lor \Diamond \Box b.
\]
Together with (\ref{Eq:Free_trick1}), this implies that $\Box f(c) = \Box \Diamond \Box b$. Moreover, we have that
\[
f(\Box c) = f (\Box \Diamond \Box a) = \Box \Diamond \Box b.
\]
Hence we conclude that $\Box f(c) = \Box \Diamond \Box b = f(\Box c)$.

Cases 2 and 3 are almost straightforward. For this reason we detail only case 2. Consider $c \in [\Box a, a]$. Since $f$ is monotone, we have that $\Box b = f(\Box a) \leq f(c) \leq f(a) = b$. In particular, this implies that $\Box f(c) = \Box b$. Moreover, we have that $f(\Box c) = f(\Box a) = 	\Box b$. Hence we conclude that $\Box f(c) = \Box b = f(\Box c)$.

We conclude that $f$ preserves $\Box$. A similar argument shows that $f$ preserves $\Diamond$ as well. Hence $f \colon \A \to \B$ is a homomorphism, as desired.
\end{proof}

It is natural to wonder whether the free two-generated positive S4-algebra is finite as well. It turns out that this is not the case, as witnessed by a small modification of the example in \cite[Figure 4]{BezGrigaml01}:

\begin{Lemma}
The free two-generated positive S4-algebra is infinite.
\end{Lemma}

\begin{proof}
It will be enough to show that there exists an infinite two-generated positive S4-algebra. To this end, consider the set of natural numbers $\omega$ equipped with the relation $\geq$. Clearly $\geq$ is a reflexive and transitive relation on $\omega$. This implies that the structure
\[
\A = \langle \mathcal{P}(\omega), \cap, \cup, \Box_{\geq}, \Diamond_{\geq}, \emptyset, \omega \rangle
\]
is indeed a positive S4-algebra.

We define recursively a formula $\varphi_{n}(x, y)$, for every $n \in \omega$, as follows:
\[
\varphi_{0} \coloneqq \Box x \text{ and }\varphi_{m+1} = \left\{\begin{array}{@{\,}ll}
\Box( x \lor  \varphi_{m}) & \text{if $m$ is odd}\\
\Box( y \lor  \varphi_{m}) & \text{if $m$ is even.}\\
\end{array} \right.
\]
Then consider the following elements of $\A$:
\[
c \coloneqq \{ 2n : n \in \omega \} \text{ and }d \coloneqq \{ 2n+1 : n \in \omega \}.
\]
We have that for every $n \in \omega$,
\[
\varphi_{n}^{\A}(c, d) = \{ m  \in \omega : m \leq n \}.
\]
Thus if $n \ne m$, then $\varphi_{n}^{\A}(c, d) \ne \varphi_{m}^{\A}(c, d)$. Hence the subalgebra $\B$ of $\A$, generated by the two-element set $\{ c, d \}$, is infinite.
\end{proof}

A variety $\class{K}$ is \textit{locally finite} if its finitely generated members are finite. Equivalently $\class{K}$ is locally finite, when its finitely generated free algebras are finite.

\begin{Corollary}\label{Cor:Non-locfinite}
$\class{PS4}$ is not locally finite.
\end{Corollary}

\section{Bottom of the subvariety lattice}

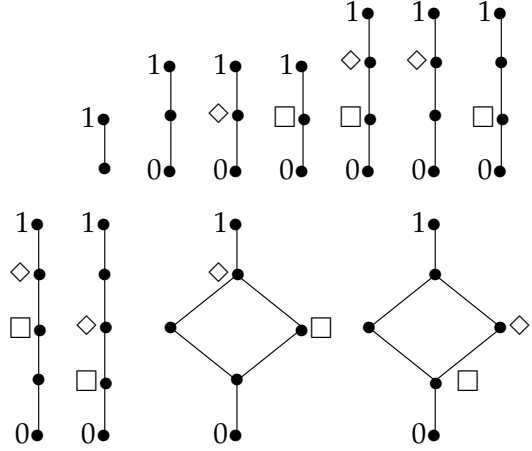
\begin{figure}
\[
\xymatrix@R=20pt @C=25pt @!0{
&&&&& *-{1 \bullet \ \ }\ar@{-}[d]  & *-{1 \bullet \ \ }\ar@{-}[d]  & *-{ 1 \bullet \ \ } \ar@{-}[d]\\
&& *-{ 1 \bullet \ \ }\ar@{-}[d]  & *-{ 1 \bullet \ \ }\ar@{-}[d]  & *-{ 1 \bullet \ \ }\ar@{-}[d]  & *-{\Diamond \bullet \ \ \?\? }\ar@{-}[d]  & *-{\Diamond \bullet \ \ \?\? }\ar@{-}[d]  & *-{\bullet}\ar@{-}[d]  \\
& *-{1 \bullet \ \ }\ar@{-}[d]  & *-{\bullet}\ar@{-}[d]  & *-{\Diamond \bullet \ \ \?\?}\ar@{-}[d]  & *-{\Box \bullet \ \ \?\? }\ar@{-}[d]  & *-{\Box \bullet \ \ \?\? }\ar@{-}[d]  &*-{\bullet}\ar@{-}[d]  &*-{\Box \bullet \ \ \?\? }\ar@{-}[d] \\
&*-{\bullet} & *-{0 \bullet \ \ } & *-{ 0 \bullet \ \ } & *-{ 0 \bullet \ \ } & *-{0 \bullet \ \ }&*-{ 0 \bullet\ \ }&*-{ 0 \bullet \ \ }\\
*-{1 \bullet \ \ }\ar@{-}[d] & *-{ 1 \bullet \ \ }\ar@{-}[d]&& *-{ 1 \bullet \ \ }\ar@{-}[d] &&& *-{ 1 \bullet \ \ }\ar@{-}[d] &\\
*-{\Diamond \bullet \ \ \?\? }\ar@{-}[d]  & *-{\bullet}\ar@{-}[d]  && *-{\Diamond \bullet \ \ \?\? }\ar@{-}[dr]\ar@{-}[dl] &&& *-{\bullet}\ar@{-}[dr]\ar@{-}[dl] &\\
*-{\Box \bullet \ \ \?\? }\ar@{-}[d]  & *-{\Diamond \bullet \ \ \?\? }\ar@{-}[d]  &  *-{\bullet}\ar@{-}[dr] &&  *-{ \ \ \?\? \bullet\Box }\ar@{-}[dl] &  *-{\bullet} \ar@{-}[dr]&&  *-{   \ \ \?\? \bullet\Diamond }\ar@{-}[dl] \\
*-{\bullet}\ar@{-}[d]  & *-{\Box \bullet \ \ \?\? }\ar@{-}[d] && *-{\bullet}\ar@{-}[d] &&& *-{\ \ \  \ \?\?  \bullet \?\? \  \Box }\ar@{-}[d] &\\
*-{0 \bullet \ \ } & *-{ 0 \bullet \ \ } && *-{ 0 \bullet \ \ } &&& *-{ 0 \bullet \ \ } &
}
\]
\caption{One-generated subdirectly irreducible algebras}\label{Fig:SIRR}
\end{figure}

In this section we will use the characterization of the free one-generated positive S4-algebra, to describe the bottom of the lattice of subvarieties of $\class{PS4}$. To this end, observe that the one-generated subdirectly irreducible positive S4-algebras coincide with the subdirectly irreducible homomorphic images of $\Fm_{\class{PS4}}(x)$. Since $\Fm_{\class{PS4}}(x)$ is described in Figure \ref{Fig:Free}, we can find these algebras by inspection. As a result we obtain that there are exactly $11$ one-generated subdirectly irreducible positive S4-algebras, which are depicted in Figure \ref{Fig:SIRR}. Reading from left to right, we will denote them respectively by
\[
\C_{2}, \D_{3}, \C_{3}^{a}, \C_{3}^{b}, \D_{4}, \C_{4}^{a}, \C_{4}^{b}, \C_{5}^{a}, \C_{5}^{b}, \C_{6}^{a} \text{ and } \C_{6}^{b}.
\]
The following observation is immediate:
\begin{Lemma}\label{Lem:MinVar}
$\VVV(\C_{2})$ is the unique minimal variety of positive S4-algebra. This variety is term-equivalent to the one of bounded distributive lattices.
\end{Lemma}

\begin{proof}
Just observe that $\C_{2}$ embeds into every non-trivial positive S4-algebra, and that $\C_{2}$ is term-equivalent to the two-element bounded distributive lattice.
\end{proof}

The next result characterizes the covers of $\VVV(\C_{2})$ in the lattice of subvarieties of $\class{PS4}$.

\begin{Theorem}\label{Thm:Covers}
$\VVV(\D_{3})$, $\VVV(\C_{3}^{a})$, $\VVV(\C_{3}^{b})$ and $\VVV(\D_{4})$ are the unique covers of $\VVV(\C_{2})$ in the lattice of subvarieties of $\class{PS4}$. Moreover, if $\class{K}$ is a subvariety of $\class{PS4}$ such that $\VVV(\C_{2}) \subsetneq \class{K}$, then $\class{K}$ includes one of these varieties.
\end{Theorem}

\begin{proof}
Applying J\'onsson's lemma to the algebras of Figure \ref{Fig:SIRR}, it is easy to see that the unique covers of $\VVV(\C_{2})$ generated by a one-generated subdirectly irreducible positive S4-algebra are $\VVV(\D_{3})$, $\VVV(\C_{3}^{a})$, $\VVV(\C_{3}^{b})$ and $\VVV(\D_{4})$. Now, consider any subvariety $\class{K}$ of $\class{PS4}$ such that $\VVV(\C_{2}) \subsetneq \class{K}$. We claim that the equational theory in one variable of $\class{K}$ differs from the one of $\C_{2}$. Suppose the contrary towards a contradiction. Then we have that
\[
\class{K} \vDash x \thickapprox \Box x \thickapprox \Diamond x.
\]
This implies that $\class{K}$ is term-equivalent to the variety of bounded distributive lattices. Keeping this in mind, it is easy to see that $\class{K} = \VVV(\C_{2})$, which contradicts the fact that $\class{K}$ properly extends $\VVV(\C_{2})$. This establishes the claim. Hence there is a one-generated subdirectly irreducible algebra $\A \ncong \C_{2}$ such that $\A \in \class{K}$. Now, we know that $\A$ is one of the algebras in Figure \ref{Fig:SIRR}. Keeping this in mind, it is easy to see that $\VVV(\A) \cap  \{ \D_{3}, \C_{3}^{a}, \C_{3}^{b}, \D_{4} \} \ne \emptyset$. Thus we conclude that $\class{K}$ contains one of the following varieties: $\VVV(\D_{3})$, $\VVV(\C_{3}^{a})$, $\VVV(\C_{3}^{b})$ and $\VVV(\D_{4})$.
\end{proof}

As it will become clear later on (Theorem \ref{Thm:StructuralCompleteness}), the variety $\VVV(\D_{4})$ is the cornerstone of structural completeness in varieties of positive K4-algebras. For this reason, we will devote to it some more attention. In particular, we will prove that $\VVV(\D_{4})$ is the variety of positive S4-algebras axiomatized by the following equations:
\begin{equation}\label{Eq:EndoEq}
\Box \Diamond x \thickapprox  \Box x\text{ and }\Diamond \Box x = \Diamond x.
\end{equation}
The proof will go through a series of lemmas.

\begin{Lemma}\label{Lem:Homomorphism}
Consider $\A \in \class{PS4}$ which satisfies (\ref{Eq:EndoEq}). The operations $\Box, \Diamond \colon A \to A$ are bounded lattice endomorphisms. Moreover, their kernels coincide and are congruences of $\A$.
\end{Lemma}

\begin{proof}
We check this for $\Box$. It is clear that $\Box$ preserve $0, 1$ and $\land$. Then consider $a, b \in A$. We have that
\begin{align*}
\Box (a \lor b) &=  \Box \Diamond ( a \lor b) = \Box ( \Diamond a \lor \Diamond b) = \Box ( \Diamond \Box a \lor \Diamond \Box b)\\
&=\Box \Diamond ( \Box a \lor \Box b) = \Box ( \Box a \lor \Box b)  = \Box a \lor \Box b.
\end{align*}
Thus we conclude that $\Box$ preserves $\lor$ and, therefore, that it is a bounded lattice endomorphism. A dual argument yields the same result for $\Diamond$.

Now observe that for every $a, b \in A$ we have that if $\Box a = \Box b$, then $\Diamond a = \Diamond \Box a = \Diamond \Box b = \Diamond b$. Thus the kernel of $\Box$ is included into the kernel of $\Diamond$. An analogous argument shows the other inclusion. Then let $\theta$ be the kernel of $\Box$ and $\Diamond$. It is clear that $\theta$ preserves $\Box$ and $\land$, since it is the kernel of $\Box$. Moreover $\theta$ preserves $\Diamond$ and $\lor$, since it is the kernel of $\Diamond$.
\end{proof}

\begin{Corollary}\label{Cor:Lattice}
Consider $\A \in \class{PS4}$ which satisfies (\ref{Eq:EndoEq}) and $a, b\in A$. If $\Box a = \Box b$, then $\textup{Cg}(a, b)$ coincides with the lattice-congruence generated by the pair $\langle a, b\rangle$.
\end{Corollary}

\begin{proof}
Let $\theta$ be the kernel of $\Box$. By Lemma \ref{Lem:Homomorphism} we know that $\theta \in \Con\A$. In particular, this implies that $\textup{Cg}(a, b) \subseteq \theta$. Let $\phi$ be the lattice-congruence generated by $\langle a, b\rangle$. We have to prove that $\phi = \textup{Cg}(a, b)$. It is clear that $\phi \subseteq \textup{Cg}(a, b)$. In order to check the other inclusion, it will be enough to show that $\phi$ preserves the modal operations. To this end consider $\langle c, d \rangle \in \phi$. Since $\phi \subseteq \textup{Cg}(a, b) \subseteq \theta$, we know that $\Box c = \Box d$ and $\Diamond c = \Diamond d$. Thus we conclude that $\langle \Box c, \Box d \rangle = \langle \Box c, \Box c\rangle \in \phi$ and, analogously, $\langle \Diamond c, \Diamond d \rangle \in \phi$.
\end{proof}

\begin{Lemma}\label{Lem:Fixed}
Consider $\A \in \class{PS4}$ subdirectly irreducible with monolith $\textup{Cg}(a, b)$ with $a < b$. If $\A$ satisfies (\ref{Eq:EndoEq}) then:
\benormal 
\item For every $c \in A$: if $b \leq c < 1$, then $c = \Diamond c$.
\item For every $c \in A$: if $0 < c \leq a$, then $c = \Box c$.
\enormal
\end{Lemma}

\begin{proof}
We prove only point 1. Suppose towards a contradiction that $c < \Diamond c$. We know that $\Box c = \Box \Diamond c$. By Corollary \ref{Cor:Lattice} we know that $\textup{Cg}(c, \Diamond c)$ coincides with the lattice-congruence generated by $\langle c, \Diamond c\rangle$. By (\ref{Eq:EDPCDL}) this means that
\[
\langle x, y \rangle \in \textup{Cg}(c, \Diamond c) \Longleftrightarrow ( x \land c = y \land c \text{ and }x \lor \Diamond c = y \lor \Diamond c).
\]
But this implies that $\langle a, b\rangle \notin \textup{Cg}(c, \Diamond c)$ and contradicts the fact that $\textup{Cg}(a, b)$ is the monolith of $\A$. We conclude that $\Diamond c = c$.
\end{proof}

In order to prove that $\VVV(\D_{4})$ is axiomatized by (\ref{Eq:EndoEq}), we need a last technical lemma.

\begin{Lemma}\label{Lem:Monoliths}
Let $\A\in \class{PS4}$ be subdirectly irreducible  satisfying (\ref{Eq:EndoEq}), and $a \in A$.
\benormal
\item If $\Box a < a$, then $\textup{Cg}(\Box a, a)$ is the monolith of $\A$.
\item If $a < \Diamond a$, then $\textup{Cg}(a, \Diamond a)$ is the monolith of $\A$.
\enormal
\end{Lemma}

\begin{proof}
Let $\textup{Cg}(b, c)$ be the monolith of $\A$ with $b < c$. We will prove only point 1.  Consider $a \in A$ such that $\Box a < a$. Applying several times point 2 of Lemma \ref{Lem:Fixed}, we obtain that
\[
a \land b = \Box ( a \land b ) = \Box a \land \Box b = b \land \Box a.
\]
Applying several times point 1 of Lemma \ref{Lem:Fixed}, we obtain that
\[
a \lor c = \Diamond ( a \lor c ) = \Diamond c \lor \Diamond a = \Diamond c \lor \Diamond \Box a = \Diamond ( c \lor \Box a) = c \lor \Box a.
\]
The two displays show that the pair $\langle \Box a, a \rangle$ belongs to the lattice-congruence generated by $\langle b, c \rangle$. In particular, this implies that $\langle \Box a, a \rangle \in \textup{Cg}(b, c)$. Since $a \ne \Box a$ and $\textup{Cg}(b, c)$ is the monolith of $\A$, we conclude that $\textup{Cg}(\Box a, a) = \textup{Cg}(b, c)$.
\end{proof}

We are finally ready to prove the desired result:

\begin{Theorem}\label{Thm:Generation}
$\VVV(\D_{4})$ is the variety of positive S4-algebras axiomatized by the equations (\ref{Eq:EndoEq}).
\end{Theorem}
\begin{proof}
Consider a subdirectly irreducible positive S4-algebra $\A$, which satisfies (\ref{Eq:EndoEq}). Moreover, suppose that $\A$ is not isomorphic to $\C_{2}$. Our goal is to show that $\A \cong \C_{4}$. First observe that for every $a \in A \smallsetminus \{ 0, 1 \}$
\begin{equation}\label{Eq:ThePoints}
\text{either }(\Box a \ne a \text{ and }a = \Diamond a) \text{ or }(\Box a = a \text{ and }a \ne \Diamond a). 
\end{equation}
To prove this, consider $a \in A \smallsetminus \{0, 1 \}$. By Corollary \ref{Cor:WellConnected} we know that either $\Box a < a$ or $a < \Diamond a$. We detail the case where $\Box a < a$, since the other one is similar. We want to prove that $a = \Diamond a$. Suppose towards a contradiction that $a < \Diamond a$. By Corollary \ref{Cor:Lattice} we know that $\textup{Cg}(\Box a, a)$ is the lattice congruence generated by the pair $\langle \Box a, a\rangle$. In particular, this implies that $\langle a, \Diamond a\rangle \notin \textup{Cg}(\Box a, a)$. But this contradicts the fact that $\textup{Cg}(a, \Diamond a)$ is the monolith of $\A$ by point 2 of Lemma \ref{Lem:Monoliths}. This establishes (\ref{Eq:ThePoints}).

Now consider the equivalence relation, defined through its blocks:
\[
\theta \coloneqq \{ \{ 0 \}, \{ \Box a : a \ne 0, 1\}, \{ \Diamond a : a \ne 0, 1 \}, \{ 1 \} \}.
\]
From (\ref{Eq:ThePoints}) it follows that the blocks of $\theta$ are pair-wise disjoint. Moreover, by (\ref{Eq:ThePoints}) we know that every element of $A$ belongs to one of the blocks. Thus $\theta$ is a well defined equivalence relation. We claim that $\theta$ is a congruence of $\A$. The fact that $\theta$ preserves the modal operations is easy to prove. We will sketch the proof of the preservation of $\Box$. Consider two \textit{different} elements $a, b$ in the same block of $\theta$. We have two cases: either $a$ and $b$ are both fixed points of $\Box$ or they are both fixed points of $\Diamond$. If $a$ and $b$ are fixed points of $\Box$, then we are done. Then consider the case where $a$ and $b$ are fixed points of $\Diamond$, i.e.,
\[
a, b \in \{ \Diamond c : c \ne 0, 1 \}.
\]
In particular, this means that $a, b \ne 0, 1$. Thus we conclude that $\Box a, \Box b \in \{ \Box c : c \ne 0, 1 \}$. We conclude that $\theta$ preserves $\Box$. A similar argument shows that $\theta$ preserves $\Diamond$.

Hence it only remains to prove that $\theta$ preserves the lattice operations. We detail the proof of the fact that $\theta$ preserves $\land$. Consider $a, b, c, d \in A$ such that $\langle a, b\rangle, \langle c, d\rangle \in \theta$. Looking at the definition of $\theta$ it is easy to see that the only non-trivial cases are the following:
\benormal
\item $a, b, c, d \in \{ \Box e : e \ne 0, 1 \}$.
\item $a, b, c, d \in \{ \Diamond e : e \ne 0, 1 \}$.
\item $a, b \in \{ \Box e : e \ne 0, 1 \}$ and $c, d \in \{ \Diamond e : e \ne 0, 1 \}$.
\enormal
1. Observe that $\Box (a \land c) = \Box a \land \Box c = a \land c$ and, similarly, $\Box (b \land d) = b \land d$. Therefore, in order to prove that $\langle a \land c, b \land d \rangle \in \theta$, it will be enough to show that $a \land c, b \land d \ne 0,1$. We already know that $a, b, c, d$ are different from $0$ and $1$. Thus $a \land c, b \land d \ne 1$. Suppose towards a contradiction that either $a \land c = 0$ or $b \land d= 0$. We assume w.l.o.g.\ that $a \land c = 0$. Applying the fact that $\Diamond$ commutes with $\land$ by Lemma \ref{Lem:Homomorphism}, we obtain that
\[
\Diamond a \land \Diamond c = \Diamond (a \land c ) = \Diamond 0 = 0.
\]
Keeping in mind that $a, c \ne 0$, this contradicts the fact that $\A$ is well-connected by Corollary \ref{Cor:WellConnected}. Thus we conclude that $a \land c \ne 0$ (and $b \land d \ne 0$). This implies that $\langle a \land c, b \land d \rangle \in \theta$.

A similar argument establishes case 2. Then we consider case 3. First observe that $a \land c$ is a fixed point of $\Box$. Suppose the contrary towards a contradiction. By Lemma \ref{Lem:Monoliths} we know that $\textup{Cg}(\Box(a \land c), a \land c)$ is the monolith of $\A$. Moreover, $\textup{Cg}(\Box(a \land c), a \land c)$ is the lattice congruence generated by the pair $\langle \Box(a \land c), a \land c\rangle$ by Corollary \ref{Cor:Lattice}. Together with the fact that $a < \Diamond a$, this implies using (\ref{Eq:EDPCDL}) that $\langle a, \Diamond a \rangle \notin \textup{Cg}(\Box(a \land c), a \land c)$. But this contradicts the fact that $\textup{Cg}(a, \Diamond a)$ is the monolith of $\A$ by Lemma \ref{Lem:Monoliths}. Thus we conclude that $a \land c$ is a fixed point of $\Box$. A similar argument shows that the same holds for $b \land d$. In order to prove that $\langle a \land c, b \land d \rangle \in \theta$ it will be enough to show that $a\land c, b \land d \ne 0, 1$. We already know that the elements $a, b, c, d$ are different from $0$ and $1$. Thus $a\land c, b \land d \ne 1$. Suppose towards a contradiction that $a \land c = 0$ or $b \land d= 0$. Assume w.l.o.g.\ that $a \land c = 0$. Keeping in mind that $\Diamond$ commutes with $\land$ by Lemma \ref{Lem:Homomorphism}, we obtain that
\[
\Diamond a \land \Diamond c = \Diamond (a \land c ) = \Diamond 0 = 0.
\]
Together with the fact that $a, c \ne 0$, this contradicts the fact that $\A$ is well-connected by Corollary \ref{Cor:WellConnected}. Hence we conclude that $\langle a \land c, b \land d \rangle \in \theta$. This establishes our claim.

Recall that the four blocks of $\theta$ are pair-wise distinct. It is not difficult to see that $\A / \theta \cong \D_{4}$. Therefore, in case $\theta = \textup{Id}_{\A}$, we are done. Suppose towards a contradiction that $\theta \ne \textup{Id}_{\A}$. By Lemma \ref{Lem:Homomorphism} we know that the kernel $\phi$ of $\Box$ coincides with the kernel of $\Diamond$ and is a congruence of $\A$. Since $\A$ is subdirectly irreducible and not isomorphic to $\C_{2}$, we know that either $\Box^{\A}$ or $\Diamond^{\A}$ is not the identity relation. This implies that $\phi \ne \textup{Id}_{\A}$. We will show that $\theta \cap \phi = \textup{Id}_{\A}$, contradicting the fact that the identity relation is meet-irreducible in $\Con \A$. Consider two different elements $a, b \in A$ such that $\langle a, b\rangle \in \phi$. This means that $\Box a = \Box b$ and $\Diamond a = \Diamond b$. Together with the fact that $a \ne b$, this implies that $a, b \ne 0, 1$. By (\ref{Eq:ThePoints}) we have two cases: either $\Box a = a$ or $a = \Diamond a$. We detail the case where $\Box a = a$. We have that $b \ne a = \Box a = \Box b$. By (\ref{Eq:ThePoints}) this means that $b = \Diamond b$. Thus $a$ and $b$ belongs to two different blocks of $\theta$, that is, $\langle a, b \rangle \notin \theta$. The case where $a = \Diamond a$ is handled similarly. Hence we obtain that $\phi \cap \theta = \textup{Id}_{\A}$, contradicting the fact that $\A$ is subdirectly irreducible. We conclude that $\A \cong \D_{4}$.

We have shown that $\C_{2}$ and $\D_{4}$ are the unique subdirectly irreducible members of the variety of positive S4-algebras axiomatized by (\ref{Eq:EndoEq}). We conclude that (\ref{Eq:EndoEq}) axiomatizes $\VVV(\D_{4})$.
\end{proof}

\section{Some splittings}

Let $\class{K}$ be a variety. A subdirectly irreducible algebra $\A \in \class{K}$ is a \textit{splitting algebra} in $\class{K}$ if there is a largest subvariety of $\class{K}$ excluding $\A$. In this section we will prove some basic results on splitting algebras that will be useful in the sequel. The following lemma is taken from \cite{McKe72} (see also \cite{Day1975}):
\begin{Lemma}[McKenzie]\label{Lem:McKenzie}
If a congruence distributive variety is generated by its finite members, then its splitting algebras are finite.
\end{Lemma}

The well-known fact that the variety of K4-algebras is generated by its finite members, together with  the fact that $\class{PKA}$ is the class of subreducts of K4-algebras (Theorem \ref{Thm:Subreducts}), implies that $\class{PKA}$ is generated by its finite members. Then its splitting algebras are finite by Lemma \ref{Lem:McKenzie}. On the other hand, we don't know which finite subdirectly irreducible algebras are splitting in $\class{PKA}$. A similar argument shows that splitting algebras in $\class{PS4}$ are finite. However, not necessarily all finite subdirectly positive S4-algebras are splitting in $\class{PS4}$.

\begin{Lemma}\label{Lem:Splitting}
\
\benormal
\item The largest subvariety of $\class{PS4}$ excluding $\C_{3}^{a}$ is axiomatized by $\Diamond \Box \Diamond x \thickapprox \Diamond x$.
\item The largest subvariety of $\class{PS4}$ excluding $\C_{3}^{b}$ is axiomatized by $\Box \Diamond \Box x \thickapprox \Box x$.
\enormal
\end{Lemma}

\begin{proof}
We will detail only the proof of 1, since the proof of 2 is similar. Let $\class{K}$ be the subvariety of $\class{PS4}$ axiomatized by $\Diamond \Box \Diamond x \thickapprox \Diamond x$. Observe that $\C_{3}^{a} \notin \class{K}$. Then suppose towards a contradiction that there is a variety $\class{W}$ such that $\C_{3}^{a} \notin\class{W}$ and $\class{W} \nsubseteq \class{K}$. This means that $\class{W} \nvDash \Diamond \Box \Diamond x \thickapprox \Diamond x$. Then there is a one-generated subdirectly irreducible algebra $\A \in \class{W}$ such that $\A \nvDash \Diamond \Box \Diamond x \thickapprox \Diamond x$. Taking a look at Figure \ref{Fig:SIRR} one sees that either $\A = \C_{3}^{a}$ or $\A = \C_{4}^{a}$. In both cases $\C_{3}^{a} \in \HHH(\A) \subseteq \class{K}$, which is false.
\end{proof}

We conclude this section by showing that $\D_{3}$ is a splitting algebra in $\class{PKA}$ (and not only in $\class{PS4}$).

\begin{Lemma}\label{Lem:Splitting2}
  $\D_{3}$ is a splitting algebras in $\class{PK4}$ and the largest subvariety of $\class{PK4}$ excluding $\D_{3}$ is axiomatized by $\Diamond x \land \Box \Diamond x \leq x \lor \Box x \lor \Diamond \Box x$.
\end{Lemma}

\begin{proof}
It is clear that a variety of positive K4-algebras satisfying $\Diamond x \land \Box \Diamond x \leq x \lor \Box x \lor \Diamond \Box x$ excludes $\D_{3}$. To prove the converse we reason towards a contradiction. Suppose that there is a variety of K4-algebras $\class{K}$, which excludes $\D_{3}$ in which the equation $\Diamond x \land \Box \Diamond x \leq x \lor \Box x \lor \Diamond \Box x$ fails. Then there is an algebra $\A \in \class{K}$ and $a \in A$ such that $\Diamond a \land \Box \Diamond a \nleq a \lor \Box a \lor \Diamond \Box a$. Consider the embedding $\kappa \colon \A \to \mathcal{M}(\A)$ defined in the proof of Theorem \ref{Thm:Subreducts}. Let $\theta$ be the congruence of $\mathcal{M}(\A)$ generated by $\langle \Box \kappa(a), 0 \rangle$ and $\langle \Diamond \kappa(a), 1 \rangle$.

We claim that $\langle \kappa(a), 1 \rangle \in \theta$. Suppose the contrary towards a contradiction. Then observe that also $\langle \kappa(a), 0 \rangle \notin \theta$, since
\[
\Diamond \kappa(a) \equiv_{\theta} 1 \not\equiv_{\theta} 0 = \Diamond 0.
\]
Then let $\phi \coloneqq \{ \langle b, c \rangle \in A^{2} : \langle \kappa(b), \kappa(c)\rangle \in \theta\}$. We have that the subalgebra of $\A / \phi$ generated by $a / \phi$ is isomorphic to $\D_{3}$. But this contradicts the fact that $\class{K}$ excludes $\D_{3}$, thus establishing the claim.

Our claim, together with (\ref{Eq:EDPCK4}) implies that
\[
(\Box \kappa(a) \leftrightarrow 0) \land \Box (\Box \kappa(a) \leftrightarrow 0) \land (\Diamond \kappa(a) \leftrightarrow 1)\land \Box(\Diamond \kappa(a) \leftrightarrow 1) \subseteq \kappa(a) \leftrightarrow 1.
\]
It is easy to see that this can be reduced to the following expression:
\begin{equation}\label{Eq:Filtters}
\lnot \Box \kappa(a) \cap \lnot \Diamond \Box \kappa(a) \cap \Diamond \kappa(a) \cap \Box  \Diamond \kappa(a) \subseteq \kappa(a).
\end{equation}
Recall that $\Diamond a \land \Box \Diamond a \nleq a \lor \Box a \lor \Diamond \Box a$ and, therefore, that
\[
\Diamond \kappa(a) \cap \Box \Diamond \kappa(a) \nsubseteq \kappa(a) \cup \Box \kappa(a) \cup \Diamond \Box \kappa(a).
\]
 By the Prime Filter Theorem, there is an ultrafilter $F \in \textup{Pr}(\mathcal{M}(\A))$ such that $\Diamond \kappa(a) \cap \Box \Diamond \kappa(a) \in F$ and $\kappa(a), \Box \kappa(a), \Diamond \Box \kappa(a) \notin F$. Since $F$ is an ultrafilter, this contradicts (\ref{Eq:Filtters}).
\end{proof}

\section{Varieties of height $\leq 4$}

In this section we characterize all varieties of positive S4-algebras that have height $\leq 4$ in the lattice of subvarieties of $\class{PS4}$. We  begin by finding the covers of $\VVV(\D_{3})$, $\VVV(\C_{3}^{a})$, $\VVV(\C_{3}^{b})$ and $\VVV(\D_{4})$ in the lattice of subvarieties of $\class{PS4}$.

\begin{Lemma}\label{Lem:Cov1}
$\VVV(\D_{3}, \D_{4})$, $\VVV(\C_{3}^{a}, \D_{4})$ and $\VVV(\C_{3}^{b}, \D_{4})$ are the unique covers of $\VVV(\D_{4})$ in the lattice of subvarieties of $\mathsf{PS4}$. Moreover, every variety $\class{K}$ of positive S4-algebras such that $\VVV(\D_{4}) \subsetneq \class{K}$ contains one of these three covers.
\end{Lemma}

\begin{proof}
By J\'onsson's lemma we obtain that $\VVV(\D_{3}, \D_{4})$, $\VVV(\C_{3}^{a}, \D_{4})$ and $\VVV(\C_{3}^{b}, \D_{4})$ are covers of $\VVV(\D_{4})$. Suppose towards a contradiction that there is a variety $\class{K}$ of positive S4-algebras such that $\VVV(\D_{4}) \subsetneq \class{K}$, which does not contain any of the varieties $\VVV(\D_{3}, \D_{4})$, $\VVV(\C_{3}^{a}, \D_{4})$ and $\VVV(\C_{3}^{b}, \D_{4})$. Then clearly $\D_{3}, \C_{3}^{a}, \C_{3}^{b} \notin \class{K}$. By Lemmas \ref{Lem:Splitting} and \ref{Lem:Splitting2} we know that $\class{K}$ satisfies the following equations:
\begin{align}
\Diamond \Box \Diamond x &\thickapprox \Diamond x \label{Eq:SPeq1}\\
\Box \Diamond \Box x &\thickapprox \Box x\label{Eq:SPeq2}\\
\Diamond x \land \Box \Diamond x &\leq x \lor \Box x \lor \Diamond \Box x\label{Eq:SPeq3}
\end{align}
Observe that in $\class{PS4}$ the inequality (\ref{Eq:SPeq3}) can be simplified as follows:
\begin{equation}\label{Eq:SPeq4}
\Box \Diamond x \leq x \lor \Diamond \Box x.
\end{equation}
It is easy to see that $\class{PS4}\vDash \Box (x \lor \Diamond \Box x) \thickapprox \Box \Diamond \Box x$. Together with (\ref{Eq:SPeq4}), this implies that
\[
\class{K} \vDash \Box \Diamond x \thickapprox \Box (\Box \Diamond x) \leq \Box (x \lor \Diamond \Box x)  \thickapprox \Box \Diamond \Box x.
\]
Together with (\ref{Eq:SPeq2}) this implies that $\class{K}\vDash \Box \Diamond x \leq \Box x$. Since $\class{PS4}\vDash \Box x \leq \Box \Diamond x$, we conclude that $\class{K}\vDash \Box \Diamond x \thickapprox \Box x$. A similar argument shows that $\class{K}\vDash \Diamond \Box x \thickapprox \Diamond x$. Hence $\class{K}$ satisfies the equations (\ref{Eq:EndoEq}). By Theorem \ref{Thm:Generation} we conclude that $\class{K} \subseteq \VVV(\D_{4})$. This contradicts the fact that $\VVV(\D_{4}) \subsetneq\class{K}$.
\end{proof}

\begin{figure}
\[
\xymatrix@R=20pt @C=25pt @!0{
&&&&&&& *-{1 \bullet \ \ }\ar@{-}[d] &&*-{1 \bullet \ \ }\ar@{-}[d]\\
& *-{1  \ \bullet \ \ \ }\ar@{-}[dr]\ar@{-}[dl] &&& *-{1 \ \bullet \ \ \ }\ar@{-}[dr]\ar@{-}[dl] &&& *-{\bullet}\ar@{-}[dl]\ar@{-}[dr]&&*-{\bullet  }\ar@{-}[d]\\
*-{\bullet  }\ar@{-}[dr]&& *-{ \bullet }\ar@{-}[dl] & *-{\bullet } \ar@{-}[dr]&& *-{\bullet  } \ar@{-}[dl]& *-{\bullet }\ar@{-}[dr] && *-{ \bullet  } \ar@{-}[dl]& *-{\bullet  }  \ar@{-}[d] \\
& *-{0 \ \bullet \ \ \  } &&& *-{ \bullet   } \ar@{-}[d]&&& *-{0 \ \bullet \ \ \ } &&*-{0 \ \bullet \ \ \  } \\
&&&&  *-{0 \ \bullet \ \ \  } &&&  &&
}
\]
\caption{Some subdirectly irreducible positive S4-algebras}\label{Fig:SIRR2}
\end{figure}
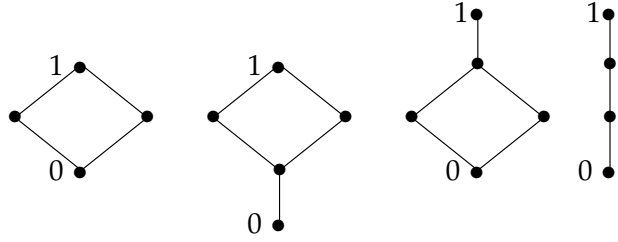

In order to describe the covers of $\VVV(\D_{3})$, $\VVV(\C_{3}^{a})$ and $\VVV(\C_{3}^{b})$ we need to introduce the new positive S4-algebras depicted in Figure \ref{Fig:SIRR2}. Reading from left to right, we denote them by $\A_{4}$, $\D_{5}^{a}$, $\D_{5}^{b}$ and $\B_{4}$.
The key observation  to understand the covers of $\VVV(\D_{3})$ is the following:

\begin{Lemma}\label{Lem.SameTheory}
Let $\class{K}$ be a variety of positive S4-algebras, whose equational theory in one variable coincides with the one of $\VVV(\D_{3})$. If $\A \in \class{K}_{si} \smallsetminus \{ \C_{2}, \D_{3}\}$, then
\[
\SSS(\A) \cap \{ \A_{4}, \B_{4} \} \ne \emptyset.
\]
\end{Lemma}

\begin{proof}
Consider $\A \in \class{K}_{si} \smallsetminus \{ \C_{2}, \D_{3}\}$. Observe that
\begin{equation}\label{Eq:NoFixedPoints}
\Box a = 0 \text{ and }\Diamond a = 1 \text{ for every }a \in A \smallsetminus \{ 0, 1 \}.
\end{equation}
To prove this consider $a \in A \smallsetminus \{ 0, 1 \}$. Observe that $\D_{3} \vDash \Box x \thickapprox \Diamond \Box x$. Since $\class{K}$ and $\D_{3}$ have the same equational theory in one variable, we conclude that $\Box a$ is a fixed point of $\Diamond$. Thus $\Box a$ is a fixed point of both modal operations. By Corollary \ref{Cor:WellConnected} we conclude that $\Box a \in \{ 0, 1 \}$. Since $a < 1$, we conclude that $\Box a = 0$. A similar argument shows that $\Diamond a = 1$. This establishes (\ref{Eq:NoFixedPoints}).

We have three cases: either $1$ is not join-prime, or $0$ is not meet-prime or none of the previous ones. First we consider the case where $1$ is not join-prime. Then there are $a, b < 1$ such that $a \lor b = 1$. From (\ref{Eq:NoFixedPoints}) it follows that the elements $\{ 0, a \land b, a, b, 1 \}$ form the universe of a subalgebra $\B$ of $\A$ such that $\B \cong \A_{4}$ or $\B \cong \D_{5}^{a}$. Observe that in both cases $\SSS(\B) \cap \{ \A_{4}, \B_{4} \} \ne \emptyset$. The case where $0$ is not meet-prime is handled similarly, showing that either $\A_{5} \in \SSS(\A)$ or $\D_{5}^{b} \in \SSS(\A)$.

It only remains to consider the case where $0$ is meet-prime and $1$ is join-prime. Together with the fact that $\A$ is different from $\C_{2}$ and $\D_{3}$, this implies that it contains a four-element chain $0 < a < b < 1$. But this is the universe of a subalgebra of $\A$ isomorphic to $\B_{4}$. Hence we conclude that $\B_{4} \in \SSS(\A)$ as desired.
%
\end{proof}

As a consequence of the above result we obtain a characterization of the covers of $\VVV(\D_{3})$.

\begin{Corollary}\label{Cor:Cov2}
$\VVV(\D_{4}, \D_{3})$, $\VVV(\C_{3}^{a}, \D_{3})$, $\VVV(\C_{3}^{b}, \D_{3})$, $\VVV(\A_{4})$, $\VVV(\B_{4})$ are the unique covers of $\VVV(\D_{3})$ in the lattice of subvarieties of $\class{PS4}$.
\end{Corollary}

\begin{proof}
From J\'onsson's lemma it follows that the varieties listed in the statement are covers of $\VVV(\D_{3})$. Consider a subvariety $\class{K}$ that covers $\VVV(\D_{3})$ in the lattice of subvarieties of $\class{PS4}$. We have two cases: either the equational theory in one variable of $\class{K}$ is the same as the one of $\VVV(\D_{3})$ or not. First consider the case where it is the same. Since $\class{K}$ covers $\VVV(\D_{3})$ there is $\A \in \class{K}_{si} \smallsetminus \{ \C_{2}, \D_{3} \}$. By Lemma \ref{Lem.SameTheory} we obtain that $\class{K} \cap \{ \A_{4}, \B_{4}\} \ne \emptyset$. Since  $\VVV(\A_{4})$ and $\VVV(\B_{4})$ are covers of $\VVV(\D_{3})$, we conclude that either $\class{K} = \VVV(\A_{4})$ or $\class{K} = \VVV(\B_{4})$. Then consider the case where the equational theory of $\class{K}$ in one variable is different from the one of $\VVV(\D_{3})$. This means that $\class{K}$ contains a one-generated subdirectly irreducible algebra $\A$ different from $\C_{2}$ and $\D_{3}$. Together with the fact that $\class{K}$ is a cover of $\VVV(\D_{3})$, this implies that $\class{K} = \VVV(\A, \D_{3})$. Now, we know that $\A$ is one of the algebras depicted in Figure \ref{Fig:SIRR}. It is easy to see that the unique algebras $\A$ in Figure \ref{Fig:SIRR} such that $\VVV(\A, \D_{3})$ is a cover of $\VVV(\D_{3})$ are $\D_{4}$, $\C_{3}^{a}$ and $\C_{3}^{b}$.
\end{proof}

It only remains to characterize the covers of $\VVV(\C_{4}^{a})$ and $\VVV(\C_{4}^{b})$. Since the algebras $\C_{4}^{a}$ and $\C_{4}^{b}$ are symmetric, we will detail the technical work only for $\VVV(\C_{4}^{a})$, since the other one is completely analogous.

\begin{Lemma}\label{Lem:SameTheory2}
$\VVV(\C_{3}^{a})$ is axiomatized w.r.t.\ $\class{PS4}$ by the equational theory of $\C_{3}^{a}$ in one variable.
\end{Lemma}

\begin{proof}
Let $\class{K}$ be a variety of positive S4-algebra, whose equational theory in one variable coincides with the one of $\C_{3}^{a}$. Our goal is to prove that $\class{K} = \VVV(\C_{3}^{a})$. First observe that $\class{K}$ properly extends $\VVV(\C_{2})$, since $\C_{3}^{a} \nvDash x \thickapprox \Box x$. Then there is $\A \in \class{K}_{si} \smallsetminus \{ \C_{2} \}$.

We will prove that any such $\A$ is isomorphic to $\C_{3}^{a}$. First observe that for every $a \in A$,
\begin{equation}\label{Eq:StrongEquations}
\Box a = 0 \text{ and }\Diamond a = a\text{ for every }a \in A \smallsetminus \{ 0, 1 \}.
\end{equation}
To prove this, consider $a \in A \smallsetminus \{ 0, 1 \}$. Observe that $\C_{3}^{a} \vDash \Diamond x \thickapprox x, \Diamond \Box x \thickapprox \Box x$. Together with the fact that the equational theory of $\class{K}$ in one variable is the same as the one of $\VVV(\C_{3}^{a})$, this implies that $\Diamond a = a$ and that $\Box a$ is a fixed point of both modal operations. From Corollary \ref{Cor:WellConnected} it follows that $\Box a \in \{ 0, 1 \}$. Together with the fact that $a < 1$, this implies that $\Box a = 0$. This establishes (\ref{Eq:StrongEquations}).

We have two cases: either $1$ is join-prime or not. If $1$ is not join-prime, then there are $a, b < 1$ such that $a \lor b = 1$. Together with  (\ref{Eq:StrongEquations}), this implies that
\[
1 = \Box (a \lor b) \leq \Box a \lor \Diamond b = 0 \lor b = b
\]
which contradicts the fact that $b < 1$. Hence we conclude that $1$ is join-prime.

Then consider $a < b$ such that $\textup{Cg}(a, b)$ is the monolith of $\A$. We have that $b < 1$. To prove this, observe that if $b=1$, then $\langle 0, 1\rangle = \langle \Box a, \Box b \rangle \in \textup{Cg}(a, b)$. This means that $\A$ is simple. Lemma \ref{Lem:Simple}, together with the fact that $\A \vDash \Diamond x \thickapprox x$, implies that $A = \{ 0, 1 \}$. But this implies that $\A = \C_{2}$, which is false by assumption. This establishes $b < 1$. Now, observe that, since $1$ is join-prime, $A \smallsetminus \{ 1 \}$ form the universe of a sublattice $\boldsymbol{L}$ of $\langle A, \land, \lor \rangle$. We know that $L$ contains at least the two elements $a < b$. We claim that $L = \{ a, b \}$. Suppose the contrary towards a contradiction. Then $\boldsymbol{L}$ is a lattice of more than three elements. In particular, this implies that there is a congruence $\theta$ of $\boldsymbol{L}$ different from the identity relation and such that $\langle a, b\rangle \notin \theta$. Consider the equivalence relation $\phi$ on $A$ whose blocks are the following:
\[
 \{ 1 \} \text{ and } c/ \theta \text{ for every } c \in L.
\]
The fact that $1$ is join-prime and (\ref{Eq:StrongEquations}) imply that $\phi$ is congruence of $\A$. Moreover, $\phi$ is different from the identity relation and $\langle a, b\rangle \notin \phi$. This contradicts the fact that $\textup{Cg}(a, b)$ is the monolith of $\A$, thus establishing the claim. Hence $A = L \cup \{ 1 \} = \{ a, b, 1 \}$. Hence $\A$ is the three-element chain, whose modal operations are determined by (\ref{Eq:StrongEquations}). We conclude that $\A \cong \C_{3}^{a}$.
\end{proof}

As a consequence of the above result we obtain a characterization of the covers of $\VVV(\C_{3}^{a})$ and $\VVV(\C_{3}^{b})$.

\begin{Corollary}\label{Cor:Cov3} \
\benormal
\item $\VVV(\C_{3}^{b}, \C_{3}^{a})$, $\VVV(\D_{3}, \C_{3}^{a})$, $\VVV(\D_{4}, \C_{3}^{a})$ and $\VVV(\C_{4}^{a})$ are the unique covers of $\VVV(\C_{3}^{a})$ in the lattice of subvarieties of $\class{PS4}$.
\item $\VVV(\C_{3}^{a}, \C_{3}^{b})$, $\VVV(\D_{3}, \C_{3}^{b})$, $\VVV(\D_{4}, \C_{3}^{b})$ and $\VVV(\C_{4}^{b})$ are the unique covers of $\VVV(\C_{3}^{b})$ in the lattice of subvarieties of $\class{PS4}$.
\enormal
\end{Corollary}

\begin{proof}
We detail the proof of $1$, since the proof of $2$ is analogous. From J\'onsson's lemma it follows that the varieties listed in the statement are covers of $\VVV(\C_{3}^{a})$. Consider a subvariety $\class{K}$ that covers $\VVV(\C_{3}^{a})$ in the lattice of subvarieties of $\class{PS4}$. By Lemma \ref{Lem:SameTheory2} we know that the equational theory of $\class{K}$ in one variable differs from the one of $\VVV(\C_{3}^{a})$. This means that $\class{K}$ contains a one-generated subdirectly irreducible algebra $\A$ different from $\C_{2}$ and $\C_{3}^{a}$. Together with the fact that $\class{K}$ is a cover of $\VVV(\C_{3}^{a})$, this implies that $\class{K} = \VVV(\A, \C_{3}^{a})$. Now, we know that $\A$ is one of the algebras depicted in Figure \ref{Fig:SIRR}. It is easy to see that the unique algebras $\A$ in Figure \ref{Fig:SIRR} such that $\VVV(\A, \D_{3})$ is a cover of $\VVV(\C_{3}^{a})$ are $\D_{4}$, $\D_{3}$, $\C_{3}^{b}$ and $\C_{4}^{a}$.
\end{proof}

Figure \ref{Fig:Varieties} describes pictorially the bottom part of the lattice of varieties of positive S4-algebras. More precisely, we have the following:

\begin{Theorem}\label{Thm:Height4}
Figure \ref{Fig:Varieties} represents all varieties of positive S4-algebras, which have height $\leq 4$ in the lattice of subvarieties of $\class{PS4}$.
\end{Theorem}

\begin{proof}
This is a consequence of Theorem \ref{Thm:Covers}, Lemmas \ref{Lem:MinVar} and \ref{Lem:Cov1}, and Corollaries \ref{Cor:Cov2} and \ref{Cor:Cov3}. 
\end{proof}

\section{Structural completeness}

We conclude this study of varieties of positive modal algebras with some application to structural completeness. To this end, some observations are in order.\ In \cite{Ryb95} Rybakov provided a full characterization of HSC varieties of K4-algebras. It turns out that they are all finitely axiomatizable and, therefore, that there are only countably many of them. Moreover, he showed that the map that associates to an intermediate logic its biggest modal companion preserves HSC \cite{MakRyb74,Ryb95}. Keeping in mind that every variety of G\"odel algebras is HSC \cite{DzWr73}, this implies that there are infinitely many HSC varieties of S4-algebras.\footnote{G\"odel algebra are simply the Heyting algebras, which are subdirect products of chains.} These considerations show that there are infinitely many (but not uncountably many) HSC varieties of S4 or, equivalently, K4-algebras. 

It is interesting to compare this situation with the one occurring in the setting of positive K4-algebras, where we will show that there are only three non-trivial SC varieties (Theorem \ref{Thm:StructuralCompleteness++}), that is SC almost never occur. Our first goal will be to prove that, in the study of SC in positive K4-algebras, we can restrict harmlessly to positive S4-algebras (Corollary \ref{Cor:K4>IA}). To this end, we will rely on the following observation:

\begin{Lemma}\label{Lem:OtherMinimal}
$\VVV(\B_{2})$ is the subvariety of $\class{PMA}$ axiomatized by $\Box x \thickapprox 1$ and $\Diamond x \thickapprox 0$.
\end{Lemma}

\begin{proof}
It is clear that the equations $\Box x \thickapprox 1$ and $\Diamond x \thickapprox 0$ holds in $\B_{2}$ and, therefore, in $\VVV(\B_{2})$. Conversely, let $\A$ be a positive modal algebra in which these equations holds. Then consider any pair of different elements $a, b \in A$. We can assume w.l.o.g.\ that $a \nleq b$. By the Prime Filter Theorem, there is $F\in \textup{Pr}(\A)$ such that $a \in F$ and $b \notin F$. Consider the map $f \colon \A \to \B_{2}$ that sends to $1 \in B_{2}$ exactly the elements which belong to $F$. Clearly $f$ is a homomorphism of bounded lattices. Moreover, since the equations $\Box x \thickapprox 1$ and $\Diamond x \thickapprox 0$ hold in $\A$, we conclude that $f$ preserves also the modal operations. Thus every pair of different elements of $\A$ is separated by a homomorphism onto $\B_{2}$. This means that $\A \in \PSD(\B_{2})$ and, therefore, that $\A \in \VVV(\B_{2})$.
\end{proof}

The next result highlights the importance of the algebra $\B_{2}$ for structural completeness in positive modal algebras.

\begin{Theorem}\label{Thm:StructuralCompletenessPMA}
Let $\class{K}$ be a SC variety of positive modal algebras. Either $\class{K} = \VVV(\B_{2})$ or there are $n, m \geq 1$ such that
\[
\class{K} \vDash \Box x \land \dots \land \Box^{n}x \leq x \text{ and }\class{K} \vDash x \leq \Diamond x \lor \dots \lor \Diamond^{m}x.
\]
\end{Theorem}
\begin{proof}
Let $\class{K}$ be a SC variety of positive modal algebras. If $\class{K}$ is trivial, then we are done. Then we can assume that $\class{K}$ is non-trivial. By Lemma \ref{Cor:PSC_simple} we know that the free $0$-generated algebra $\A$ over $\class{K}$ is simple. We claim that either $\A = \B_{2}$ or $\A=\C_{2}$. To prove this, assume that $\A \ne \B_{2}$. Then suppose towards a contradiction that $\A \ne \C_{2}$. We have that either $0 <\Box 0$ or $\Diamond 1 < 1$. We detail the case where $0 < \Box 0$. Observe that if $\Box 0=1$, then
\[
\Diamond 1 = \Diamond 1 \land \Box 0 \leq \Diamond (0 \land 1)= \Diamond 0=0.
\]
But this implies that $\A=\B_{2}$, which is false. Thus we have that $0 < \Box 0 < 1$. Consider the embedding $\kappa \colon \A \to \mathcal{M}(\A)$ described in Theorem \ref{Thm:Subreducts}. Then let $F$ be the filter of $\mathcal{M}(\A)$ generated by $\kappa(\Box 0)$. Observe that for every $b \in F$, we have that $\Box 0 \leq \Box b$ by the monotonicity of $\Box$. Thus we conclude that $F \in \textup{Op}(\mathcal{M}(\A))$. Consider the congruence $\phi \coloneqq \{ \langle b, c \rangle \in A^{2} : \langle \kappa(b), \kappa(c)\rangle \in \theta_{F} \}$ of $\A$. We have that
\[
\langle 0, \Box 0\rangle \notin \phi \text{ and }\langle 1, \Box 0\rangle \in \phi.
\]
Thus we conclude that $\phi$ is different both from the identity and from the total relation, contradicting the fact that $\A$ is simple. The case where $\Diamond 1 < 1$ is handled similarly, with the only difference that one relies on the correspondence between the congruences of $\mathcal{M}(\A)$ and its ideals closed under $\Diamond$. This establishes our claim.

Now, consider the case where $\A = \B_{2}$. Clearly in $\class{K}$ the equations $\Box x \thickapprox 1$ and $\Diamond x \thickapprox 0$ hold. By Lemma \ref{Lem:OtherMinimal} we conclude that $\class{K} = \VVV(\B_{2})$. Then consider the case where $\A \ne \B_{2}$. From the claim it follows that $\A = \C_{2}$. Suppose towards a contradiction that either there is no $n \geq 1$ such that $\class{K} \vDash \Box x \land \dots \land \Box^{n}x \leq x$ or there is no $n \geq 1$ such that $\class{K} \vDash x \leq \Diamond x \lor \dots \lor \Diamond^{n}x$. We will detail only the first case, since the second one is analogous.

Consider the free $1$-generated algebra $\Fm_{\class{K}}(x)$ over $\class{K}$. Let $\kappa \colon \Fm_{\class{K}}(x) \to \mathcal{M}(\Fm_{\class{K}}(x))$ be the embedding defined in the proof of Theorem \ref{Thm:Subreducts}. Consider the set
\[
F\coloneqq \{ a \in \mathcal{M}(\Fm_{\class{K}}(x)) : \kappa ( \llbracket \Box x \land \dots \land \Box^{n} x \rrbracket)\leq a \text{ for some }n \geq 1\}.
\]
It is easy to see that $F$ is a filter of $\mathcal{M}(\Fm_{\class{K}}(x))$. Moreover, consider an arbitrary $a \in F$. Then there is $n \geq 1$ such that $\kappa ( \llbracket \Box x \land \dots \land \Box^{n} x \rrbracket)\leq a$. We have that:
\begin{align*}
\kappa (\llbracket \Box^{1} x \land \dots \land \Box^{n+1} x \rrbracket)&\leq \kappa (\llbracket \Box^{2} x \land \dots \land \Box^{n+1} x \rrbracket)\\
&=\kappa (\Box \llbracket \Box x \land \dots \land \Box^{n} x \rrbracket)\\
&=\Box\kappa ( \llbracket \Box x \land \dots \land \Box^{n} x \rrbracket)\\
&\leq \Box a.
\end{align*}
We conclude that $\Box a \in F$. This shows that $F$ is indeed an open filter. Then let $\phi \coloneqq \{ \langle b, c \rangle \in Fm_{\class{K}}(x)^{2} : \langle \kappa(b), \kappa(c)\rangle \in \theta_{F} \}$. Since $\kappa(\llbracket x \rrbracket)  \notin F$ and $\kappa(\llbracket \Box x \rrbracket) \in F$, we have that
\[
\Box (\llbracket x\rrbracket / \phi) = \llbracket 1 \rrbracket / \phi \text{ and }\llbracket x\rrbracket / \phi < \llbracket 1 \rrbracket / \phi.
\]
By Birkhoff's Subdirect Representation Theorem, there is an onto homomorphism $f \colon \Fm_{\class{K}}(x)/ \phi \to \B$ where $\B$ is subdirectly irreducible and $f(\llbracket x\rrbracket / \phi)< f( \llbracket 1 \rrbracket / \phi)$. By condition  3 of Theorem \ref{Thm:Bergman}, we know that $\B \in \SSS\PPU(\Fm_{\class{K}}(\omega))$. This means that there is an ultrapower of $\Fm_{\class{K}}(\omega)$ in which the following sentence holds:
\[
\exists x (x < 1 \text{ and }\Box x \thickapprox 1).
\]
By \L os's Theorem \cite[Theorem V.2.9]{BuSa00}, the above sentence holds in $\Fm_{\class{K}}(\omega)$. Thus there is a term $t(x_{1}, \dots, x_{n})$ such that
\begin{equation}\label{Eq:BoxEquations}
\class{K}\vDash \Box t(x_{1}, \dots, x_{n}) \thickapprox 1 \text{ and }\class{K} \nvDash t(x_{1}, \dots, x_{n}) \thickapprox 1.
\end{equation}
Then there is an algebra $\C \in \class{K}$ and $a_{1}, \dots, a_{n} \in C$ such that $t^{\C}(a_{1}, \dots, a_{n}) < 1$. Observe that all the basic  operations of $\C$ are monotone. Thus an easy induction on the construction of $t$ shows that
\[
t^{\C}(0, \dots, 0) \leq t^{\C}(a_{1}, \dots, a_{n}) < 1.
\]
Now, observe that $\B_{2}$ is a two-element subalgebra of $\C$. Hence we conclude that $t^{\C}(0, \dots, 0) = 0$. In particular, this implies that $\Box^{\C}t^{\C}(0, \dots, 0) = \Box^{\C}0 = 0$. But this contradicts (\ref{Eq:BoxEquations}). 
\end{proof}

\begin{Corollary}\label{Cor:K4>IA}
Let $\class{K}$ be a SC variety of positive K4-algebras. Either $\class{K} = \VVV(\B_{2})$ or $\class{K} \subseteq \class{PS4}$.
\end{Corollary}

\begin{proof}
Suppose that $\class{K} \ne \VVV(\B_{2})$. By Theorem \ref{Thm:StructuralCompletenessPMA} there are $n, m \geq 1$ such that 
\[
\class{K} \vDash \Box x \land \dots \land \Box^{n}x \leq x \text{ and }\class{K} \vDash x \leq \Diamond x \lor \dots \lor \Diamond^{m}x.
\]
Since $\class{K} \subseteq \class{PK4}$, it satisfy $\Box x \land \dots \land \Box^{n}x \thickapprox \Box x$ and  $\Diamond x \lor \dots \lor \Diamond^{m}x \thickapprox \Diamond x$. Hence we conclude that the equations  $\Box x \leq x \leq \Diamond x$ hold in $\class{K}$.
\end{proof}

In the light of Corollary \ref{Cor:K4>IA}, the only difficult part in understanding structural completeness in positive \textit{K4-algebras} consists in characterizing the structurally complete varieties of \textit{S4-algebras}. We begin by proving some negative results.

\begin{Lemma}\label{Lem:C3}
Every ASC variety of positive S4-algebras excludes $\C_{3}^{a}$ and $\C_{3}^{b}$.
\end{Lemma}
\begin{proof}
Suppose towards a contradiction that there is an ASC variety $\class{K} \subseteq \class{PS4}$ that contains either $\C_{3}^{a}$ or $\C_{3}^{b}$. We detail the case where $\C_{3}^{b} \in \class{K}$, since the other one is similar. By condition 1 of Theorem \ref{Thm:Bergman}, we know that $\C_{3}^{b} \times \C_{2} \in \QQQ(\Fm_{\class{K}}(\omega))$. Observe that the congruence lattice of $\C_{3}^{b} \times \C_{2}$ has exactly four elements, namely the identity, the total relation and the two congruences corresponding to the projections associated with the product. In particular, this implies that either $\C_{3}^{b} \in \QQQ(\Fm_{\class{K}}(\omega))$ or $\C_{3}^{b} \times \C_{2}$ is subdirectly irreducible relative to $\QQQ(\Fm_{\class{K}}(\omega))$. Recall that $\C_{3}^{b}$ is subdirectly irreducible in the absolute sense. Thus we conclude that
\[
\text{either }\C_{3}^{b} \in \SSS\PPU(\Fm_{\class{K}}(\omega))\text{ or }\C_{3}^{b} \times \C_{2} \in \SSS\PPU(\Fm_{\class{K}}(\omega)).
\]
Consider the case where $\C_{3}^{b} \in \SSS\PPU(\Fm_{\class{K}}(\omega))$. This means that $\C_{3}^{b}$ embeds into an ultrapower of $\Fm_{\class{K}}(\omega)$, where the following first-order sentence holds:
\[
\exists x (0 \ne x \ne 1 \textrm{ and } \Box x \thickapprox x \textrm{ and } \Diamond x \thickapprox 1).
\]
By \L os's Theorem \cite[Theorem V.2.9]{BuSa00}, the above sentence holds in $\Fm_{\class{K}}(\omega)$ too. In particular, this implies that $\C_{3}^{b} \in \SSS(\Fm_{\class{K}}(\omega))$. Similarly, one shows that if $\C_{3}^{b} \times \C_{2} \in \SSS\PPU(\Fm_{\class{K}}(\omega))$, then $\C_{3}^{b} \times \C_{2} \in \SSS(\Fm_{\class{K}}(\omega))$ as well. Hence we conclude that
\[
\text{either }\C_{3}^{b} \in \SSS(\Fm_{\class{K}}(\omega))\text{ or }\C_{3}^{b} \times \C_{2} \in \SSS(\Fm_{\class{K}}(\omega)).
\]
In both cases, this implies that there is $\llbracket \varphi(x_{1}, \dots, x_{n}) \rrbracket \in  Fm_{\class{K}}(\omega)$ such that $\llbracket \varphi \rrbracket < \llbracket 1 \rrbracket$ and $\Diamond \llbracket \varphi \rrbracket = \llbracket 1 \rrbracket$. Now it is easy to prove by induction on terms that
\[
\llbracket \psi \rrbracket \leq \llbracket \Diamond x_{1} \lor \dots \lor \Diamond x_{n} \rrbracket
\]
for every $\llbracket \psi \rrbracket\in Fm_{\class{K}}(x_{1}, \dots, x_{n})$ such that $\llbracket \psi \rrbracket \ne \llbracket 1 \rrbracket$. In particular, this implies that $\llbracket \varphi \rrbracket \leq \llbracket \Diamond x_{1} \lor \dots \lor \Diamond x_{n} \rrbracket$ and, therefore, that
\[
\llbracket 1 \rrbracket = \Diamond \llbracket \varphi \rrbracket \leq \Diamond \llbracket \Diamond x_{1} \lor \dots \lor \Diamond x_{n} \rrbracket= \llbracket \Diamond x_{1} \lor \dots \lor \Diamond x_{n} \rrbracket.
\]
Thus we conclude that
\begin{equation}\label{Eq:DiamondToOne}
\class{K} \vDash \Diamond x_{1} \lor \dots \lor \Diamond x_{n} \thickapprox 1.
\end{equation}
Observe the equation in (\ref{Eq:DiamondToOne}) is not valid in $\C_{3}^{b}$, since $\Diamond 0 \lor\dots \lor \Diamond 0 = 0 \ne 1$. Thus we obtain a contradiction as desired.
\end{proof}

Another negative result that turns out to be very useful is the following:

\begin{Lemma}\label{Lem:D3}
Every ASC variety of $\class{PS4}$ excludes $\D_{3}$.
\end{Lemma}

\begin{proof}
Suppose the contrary towards a contradiction that $\class{K}$ is a ASC variety of $\class{PS4}$ that contains $\D_{3}$. By condition 1 of Theorem \ref{Thm:Bergman} we know that $\D_{3} \times \C_{2} \in \QQQ(\Fm_{\class{K}}(\omega))$. An argument similar to the one applied in the proof of Lemma \ref{Lem:C3} shows that there is $n \in \omega$ such that either $\D_{3}$ or $\D_{3} \times \C_{2}$ embeds into $\Fm_{\class{K}}(x_{1}, \dots, x_{n})$. In both cases, this implies that there is $\llbracket \varphi \rrbracket \in  Fm_{\class{K}}(x_{1}, \dots, x_{n})$ such that $\llbracket \varphi \rrbracket < \llbracket 1 \rrbracket$ and $\Diamond \llbracket \varphi \rrbracket = \llbracket 1 \rrbracket$. Again, as in the proof of Lemma \ref{Lem:C3}, this yields that $\class{K}\vDash \Diamond x \thickapprox 1$ which is false.
\end{proof}

We are now ready to prove the main results of this section:

\begin{Theorem}\label{Thm:StructuralCompleteness}
Let $\class{K}$ be a non-trivial variety of positive S4-algebras. The following conditions are equivalent:
\benroman
\item $\class{K}$ is actively structurally complete.
\item $\class{K}$ is structurally complete.
\item $\class{K}$ is hereditarily structurally complete.
\item $\class{K} = \VVV(\C_{2})$ or $\class{K} = \VVV(\D_{4})$.
\item $\class{K}$ satisfies the equations  $\Box \Diamond x \thickapprox \Box x$ and $\Diamond \Box x \thickapprox \Diamond x$.
\eroman
\end{Theorem}
\begin{proof}
Clearly (iii)$\Rightarrow$(ii) and (ii)$\Rightarrow$(i). Moreover, the equivalence between (iv) and (v) is the content of Theorem \ref{Thm:Generation}. 

(i)$\Rightarrow$(iv). By Lemmas \ref{Lem:C3} and \ref{Lem:D3} we know that $\class{K}$ excludes $\C_{3}^{a}, \C_{3}^{b}$ and $\D_{3}$. Together with Theorem \ref{Thm:Covers}, this implies that either $\class{K} = \VVV(\C_{2})$ or $\D_{4} \in \class{K}$. If $\class{K} = \VVV(\C_{2})$, we are done. Then consider the case where $\D_{4} \in \class{K}$. Lemma \ref{Lem:Cov1}, together with the fact that $\class{K}$ excludes  $\C_{3}^{a}, \C_{3}^{b}$ and $\D_{3}$, implies that $\class{K} = \VVV(\D_{4})$.

(iv)$\Rightarrow$(iii). By condition 4 of Theorem \ref{Thm:Bergman} We have to prove that every subquasi-variety of $\VVV(\D_{4})$ is a variety. Observe that the one-generated free algebra $\A$ in $\VVV(\D_{4})$ is the five-element chain $0 < a < b < c < 1$ where $a = \Box b = \Box c$ and $c = \Diamond a = \Diamond b$. It is easy to see that $\D_{4}$ is a retract of $\A$. Thus we conclude that $\D_{4}$ is projective in $\VVV(\D_{4})$.

Now consider a non-trivial subquasi-variety $\class{K}$ of $\VVV(\D_{4})$. Clearly $\VVV(\class{K}) = \VVV(\C_{2})$ or $\class{K} = \VVV(\D_{4})$. First consider the case where $\class{K} = \VVV(\C_{2})$. Clearly $\C_{2} \in \class{K}$, since $\class{K}$ is non-trivial. Thus we have that
\[
\VVV(\C_{2}) = \PSD(\C_{2})  \subseteq \PSD(\class{K})= \class{K}.
\]
Hence we conclude that $\class{K} = \VVV(\C_{2})$. Second consider the case where $\class{K} = \VVV(\D_{4})$. This means that $\D_{4} \in \HHH(\class{K})$. Then there is an algebra $\B \in \class{K}$ such that $\D_{4} \in \HHH(\B)$. Since $\D_{4}$ is projective in $\VVV(\D_{4})$, we conclude that $\D_{4} \in \SSS(\B) \subseteq \class{K}$. Thus we have that
\[
\VVV(\D_{4}) = \PSD(\C_{2}, \D_{4}) \subseteq \PSD(\class{K}) = \class{K}.
\]
This implies that $\class{K}= \VVV(\D_{4})$ and, therefore, we are done.
\end{proof}

As a consequence we obtain a characterization of structurally complete varieties of positive K4-algebras:

\begin{Theorem}\label{Thm:StructuralCompleteness++}
Let $\class{K}$ be a non-trivial variety of positive K4-algebras. The following conditions are equivalent:
\benroman
\item $\class{K}$ is structurally complete.
\item $\class{K}$ is hereditarily structurally complete.
\item $\class{K}= \VVV(\B_{2})$ or $\class{K} = \VVV(\C_{2})$ or $\class{K} = \VVV(\D_{4})$.
\eroman
\end{Theorem}

\begin{proof}
Part (ii)$\Rightarrow$(i) is clear. (i)$\Rightarrow$(iii). Suppose that $\class{K} \ne \VVV(\B_{2})$. Then $\class{K} \subseteq \class{PS4}$ by Corollary \ref{Cor:K4>IA}. Hence we can apply Theorem \ref{Thm:StructuralCompleteness} obtaining that $\class{K} = \VVV(\C_{2})$ or $\class{K} = \VVV(\D_{4})$.

(iii)$\Rightarrow$(ii). By Theorem \ref{Thm:Representation} we know that $\VVV(\C_{2})$ and $\VVV(\D_{4})$ are hereditarily structurally complete. Then consider the case where $\class{K} = \VVV(\B_{2})$. Observe that $\B_{2} \in \SSS(\A)$ for every non-trivial $\A \in \VVV(\B_{2})$. This implies that the only subquasi-varieties of $\VVV(\B_{2})$ are the trivial and the total ones. Hence we conclude that $\VVV(\B_{2})$ is hereditarily structurally complete.
\end{proof}

\begin{problem}
Are there ASC varieties of positive K4-algebras that are not SC?
\end{problem}

We conclude this section by characterizing passive structural completeness in varieties of positive K4-algebras.

\begin{Theorem}\label{Thm:PSC}
Let $\class{K}$ be a non-trivial variety of positive K4-algebras. The following conditions are equivalent:
\benroman
\item $\class{K}$ is passively structurally complete.
\item Either $\class{K} = \VVV(\B_{2})$ or ($\Fm_{\class{K}}(0) = \C_{2}$ and $\C_{2}$ is the unique simple member of $\class{K}$).
\item Either $\class{K} = \VVV(\B_{2})$ or ($\Fm_{\class{K}}(0) = \C_{2}$ and $\class{K}$ excludes $\D_{3}$).
\item Either $\class{K} = \VVV(\B_{2})$ or
\[
\class{K}\vDash \Diamond 1 \thickapprox 1, \Box 0 \thickapprox 0,\Diamond x \land \Box \Diamond x \leq x \lor \Box x \lor \Diamond \Box x.
\]
\eroman
\end{Theorem}
\begin{proof}
(i)$\Rightarrow$(iii). The same argument described in the proof of Theorem \ref{Thm:StructuralCompletenessPMA} shows that $\Fm_{\class{K}}(0)$ is either $\B_{2}$ or $\C_{2}$. If $\Fm_{\class{K}}(0) = \B_{2}$, then clearly $\class{K} \vDash \Box x \thickapprox 1, \Diamond x \thickapprox 0$. By Lemma \ref{Lem:OtherMinimal} this implies that $\class{K} = \VVV(\B_{2})$. Then consider the case where $\Fm_{\class{K}}(0) = \C_{2}$. We have to prove that  $\class{K}$ excludes $\D_{3}$. Suppose the contrary towards a contradiction. Then the following positive existential sentence holds in $\D_{3}$:
\[
\exists x (0 \thickapprox \Box x \textrm{ and } \Diamond x \thickapprox 1).
\]
By condition 2 of Theorem \ref{Thm:Bergman} we conclude that the above sentence holds in every non-trivial member of $\class{K}$. But this contradicts the fact that $\C_{2} \in \class{K}$, since $\class{K}$ is non-trivial.

(iii)$\Rightarrow$(ii). Suppose that $\class{K} \ne \VVV(\B_{2})$. Then consider a simple algebra $\A \in \class{K}$. By Corollary \ref{Cor:WellConnected} we know that $\Box a = 0$ and $\Diamond a = 1$ for every $a \in A \smallsetminus \{ 0, 1 \}$. This implies that if $\A$ has at least three elements, then $\D_{3} \in \SSS(\A)$. From the assumption we conclude that $\A$ has only two elements. In particular, this implies that $\A \cong \C_{2}$.

(ii)$\Rightarrow$(i). We detail the case where $\class{K} \ne \VVV(\B_{2})$, since the other one is analogous (and easier). By condition 2 of Theorem \ref{Thm:Bergman} it will be enough to show that if a positive existential sentence $\Phi$ holds in a non-trivial member of $\class{K}$, then $\Phi$ holds in every non-trivial member of $\class{K}$. To this end, consider such a sentence $\Phi$ and suppose that it holds in a non-trivial $\A \in \class{K}$. Since $\Phi$ is a positive existential sentences, it has the form
\[
\Phi = \exists x_{1}, \dots, x_{n} \Psi
\]
where $\Psi$ is a conjunction of disjunctions of equations in variables $x_{1}, \dots, x_{n}$. We know that there are elements $a_{1}, \dots, a_{n}$ such that $\A \vDash \Psi(a_{1}, \dots, a_{n})$. Then let $\B$ be the subalgebra of $\A$ generated by $a_{1}, \dots, a_{n}$. It is clear that $\Phi$ holds in $\B$ too. It is well known that every non-trivial finitely generated algebra of finite type has a simple homomorphic image \cite[pp. 153-154]{Jo70}. Together with the assumptions, this implies that $\C_{2} \in \HHH(\B_{2})$. Since positive existential sentences are preserved under homomorphic images, we conclude that $\Phi$ holds in $\C_{2}$. Since $\C_{2}$ embeds into every non-trivial member of $\class{K}$, we obtain that $\Phi$ holds in every non-trivial member of $\class{K}$ as desired.

The equivalence between (iii) and (iv) follows from Lemma \ref{Lem:Splitting2}.
\end{proof}

The next example shows that there are infinitely many varieties of positive S4-algebras both with and without PSC.

\begin{exa}\label{Ex:PSC}
We begin by proving that there are infinitely many PSC varieties of positive S4-algebras . For every natural $n \geq 1$, consider the power-set Boolean algebra $\mathcal{P}(\{ a_{1}, \dots, a_{n} \})$. We equip it with a modal operator $\Box$ defined as follows:
\begin{displaymath}
\Box X = \left\{\begin{array}{@{\,}ll}
X & \text{if $X = \{ a_{1}, \dots, a_{n}\}$}\\
\emptyset & \text{if $a_{1} \notin X$}\\
\{ a_{1} \} & \text{otherwise}\\
\end{array} \right.
\end{displaymath}
for every $X \subseteq \{ a_{1}, \dots, a_{n} \}$. It is easy to see that the resulting expansion $\A_{n}$ is an S4-algebra. We denote by $\A_{n}^{-}$ its positive reduct (with $\Diamond$ in the signature). Clearly $\A_{n}^{-} \in \class{PS4}$. Applying the correspondence between congruences and open filters in S4-algebras to the definition of $\A_{n}$, we conclude that $\A_{n}$ is subdirectly irreducible. Moreover, the fact that $\langle A, \land, \lor, 0, 1\rangle$ is a Boolean lattice implies that the congruences of $\A_{n}^{-}$ preserve the set-theoretic complement operation. Thus the congruences of $\A_{n}$ and $\A_{n}^{-}$ coincide. In particular, this implies that the identity relation in meet-irreducible in $\Con\A_{n}^{-}$ and, therefore, that $\A_{n}^{-}$ is subdirectly irreducible.

Consider $n < m$.  We know that $\A_{n}^{-}$ is subdirectly irreducible. Thus we can apply Jonsson's lemma, that on cardinality grounds yields that $\A_{m}^{-} \notin \VVV(\A_{n}^{-})$. Hence we conclude that if $n \ne m$, then $\VVV(\A_{n}^{-}) \ne \VVV(\A_{m}^{-})$. Thus there are infinitely many varieties of the form $\VVV(\A_{n}^{-})$ with $n \geq 1$. It only remains to prove that these varieties are PSC. By Theorem \ref{Thm:PSC} it will be enough to show that $\VVV(\A_{n}^{-})$ excludes $\D_{3}$ for every $n \geq 1$. By J\'onsson's lemma, this amounts to proving that $\D_{3} \notin \HHH\SSS(\A_{n}^{-})$. Suppose towards a contradiction that $\D_{3} \in \HHH\SSS(\A_{n}^{-})$, i.e., that there is a subalgebra $\B$ of $\A_{n}^{-}$ and a congruence $\theta \in \Con \B$ such that $\D_{3} \cong \B/ \theta$. Thus $\B/ \theta$ is a three element chain. Let $b/ \theta$ be the intermediate element of this chain. We know that $\langle \Box b, 0 \rangle, \langle \Diamond b, 1\rangle \in \theta$. Recall that $b \subseteq \{ a_{1}, \dots, a_{n} \}$. Moreover, observe that $\emptyset \ne b \ne \{ a_{1}, \dots, a_{n} \}$, otherwise $\B / \theta$ would be the trivial algebra. We have two cases: either $a_{1}\in b$ or $a_{1} \notin b$. First we consider the case where $a_{1} \in b$. By definition of $\A_{n}^{-}$ we have that
\[
\Diamond \Box b = \Diamond \{ a_{1} \} = \{ a_{1}, \dots, a_{n} \} = 1^{\A^{-}_{n}}.
\]
But this implies that $\langle 1, 0 \rangle = \langle \Diamond \Box b, \Diamond 0\rangle \in \theta$, contradicting the fact that $\B / \theta$ is non-trivial. Then we consider the case where $a_{1} \notin b$. Again the definition of $\A_{n}^{-}$ implies that
\[
\Box \Diamond b = \Box \{ a_{2}, \dots, a_{n} \} = \emptyset = 0^{\A^{-}_{n}}.
\]
But this implies that $\langle 0, 1 \rangle = \langle \Box \Diamond b, \Box 1\rangle \in \theta$, contradicting the fact that $\B / \theta$ is non-trivial. Hence we conclude that $\D_{3} \notin \HHH\SSS(\A_{n}^{-})$ as desired. This shows that there are infinitely many PSC varieties of positive S4-algebras.

In order to construct infinitely many varieties of positive S4-algebras, that are not PSC, we reason as follows. For every $n \geq 2$, let $\A_{n}$ be the positive S4-algebra whose universe is the $2^{n}$-element Boolean lattice and whose modal operations are defined as in condition (i) of Lemma \ref{Lem:Simple}. By Lemma \ref{Lem:Simple} we know that each $\A_{n}$ is simple. Observe that Jonsson's lemma implies, on cardinality grounds, that if $n < m$  yields that $\A_{m}^{-} \notin \VVV(\A_{n})$.  Thus there are infinitely many varieties of the form $\VVV(\A_{n}^{-})$ with $n \geq 2$. Moreover, it is clear that $\D_{3} \in \SSS(\A_{n}) \subseteq \VVV(\A_{n})$ for every $n \geq 2$. By Theorem \ref{Thm:PSC} we conclude that the varieties of the form $\VVV(\A_{n})$ are not PSC.
\qed
\end{exa}

\section*{Appendix}\label{Sec:Gentzen}

 Let $\vdash_{{\bf K}}$ be the \textit{local} consequence relation \cite{MK07c} associated with the normal modal logic ${\bf K}$. Then let $\vdash^{+}_{{\bf K}}$ be the restriction of $\vdash_{{\bf K}}$ to the set of positive modal formulas. The  logic $\vdash^{+}_{{\bf K}}$ has been studied in \cite{CeJa97,Du95} under the name of \textit{positive modal logic} (for short PML). In particular, in \cite{CeJa97} PML was axiomatized by means of the Gentzen system $\mathcal{PML}$ mentioned in the introduction. The system $\mathcal{PML}$ was formulated as a consequence relation $\vdash_{\mathcal{PML}}$ on the following set of sequents:
\[
Seq \coloneqq \{ \Gamma \rhd \varphi : \Gamma \cup \{ \varphi \} \text{ is a finite set of positive modal formulas}\}.
\]
The relation between PML and $\vdash_{\mathcal{PML}}$ is as follows: for every set of formulas $\Gamma \cup \{ \varphi \}$ we have that
\[
\Gamma \vdash^{+}_{{\bf K}} \varphi \Longleftrightarrow \text{ there is a finite }\Delta \subseteq \Gamma \text{ s.t. }\emptyset \vdash_{\mathcal{PML}} \Delta \rhd \varphi.
\]
The Gentzen system $\mathcal{PML}$ enjoys a strong relation with the variety of positive modal algebras. In order to explain this, we need to introduce some definition. Given a sequent $\Gamma\rhd \varphi$ in $Seq$, where $\Gamma = \{ \gamma_{1}, \dots, \gamma_{n} \}$, we define
\[
\btau(\Gamma\rhd \varphi) \coloneqq  \{ \gamma_{1}\land \dots \land \gamma_{n} \leq \varphi \}.
\]
Moreover, given an equation $\varphi \thickapprox \psi$ in the language of positive modal algebras, we define
\[
\brho(\varphi \thickapprox \psi) \coloneqq \{ \varphi \rhd \psi, \psi \rhd \varphi \}.
\]
The maps $\btau$ and $\brho$ can be extended to power sets by picking unions. In this way we obtain two maps
\[
\btau \colon \mathcal{P}(Seq) \longleftrightarrow \mathcal{P}(Eq) \colon \brho
\]
where $Eq$ is the set of equations in the language of positive S4-algebras. It turns out that the consequence relation $\vdash_{\mathcal{PML}}$ and the equational consequence $\vDash_{\class{PMA}}$ relative to $\class{PMA}$ are mutually interpretable in the following sense \cite{Ja02}. For every $\Sigma \cup \{ \Gamma\rhd \varphi \} \subseteq Seq$ and every $\varphi \thickapprox \psi \in Eq$ we have that:
\begin{align*}
\Sigma \vdash_{\mathcal{PML}}\Gamma\rhd \varphi &\Longleftrightarrow \btau(\Sigma) \vDash_{\class{PMA}}\btau(\Gamma \rhd \varphi) \\
\varphi\thickapprox \psi &\semeq_{\class{PMA}} \btau\brho(\varphi \thickapprox \psi).
\end{align*}
In abstract algebraic logic \cite{AAL-AIT-f,Cz01,FJa09} the relation described above is known as \textit{algebraizability}. More precisely, the Gentzen system $\mathcal{PML}$ is algebraizable in the sense of \cite{ReV95-p} (see also \cite{BP89,FJa09,Ra06}) with equivalent algebraic semantics the variety $\class{PMA}$. This induces a dual isomorphism between the lattices of \textit{varieties} of positive modal algebras and \textit{axiomatic extensions} of $\mathcal{PML}$, i.e. extensions of $\mathcal{PML}$ by means of rules of the form $\emptyset \vdash \Gamma \rhd \varphi$. In particular, the variety of postive K4-algebras correspond to the axiomatic extension $\mathcal{PML}_{4}$, obtained adding to $\mathcal{PML}$ the axioms
\[
\emptyset \vdash \brho (\Box x \leq \Box \Box x) \text{ and } \emptyset \vdash \brho (\Diamond \Diamond x \leq \Diamond x).
\]
Observe that the dual isomorphism between varieties of positive modal algebras and axiomatic extensions of $\mathcal{PML}$ preserves and reflects (active, passive) structural completeness. Thus Theorems \ref{Thm:StructuralCompleteness++} and \ref{Thm:PSC} provide a characterization of the various kinds of structural completeness in the axiomatic extensions of $\mathcal{PML}_{4}$.

\

\paragraph{\bfseries Acknowledgements.}
We wish to thank to Prof.\ Josep Maria Font for providing several useful comments on a first draft of this paper. Moreover, special thanks are due to the anonymous referees for providing additional references and useful remarks, which helped to improve significantly the presentation. This research was supported by the joint project of Austrian Science Fund (FWF) I$1897$-N$25$ and Czech Science Foundation (GACR) GF$15$-$34650$L, and by project CZ.$02$.$2$.$69$/$0$.$0$/$0$.$0$/$17$\_$050$/$0008361$, OPVVV M\v{S}MT, MSCA-IF Lidsk\'{e} zdroje v teoretick\'{e} informatice.

\begin{landscape}
\begin{figure}
\[
\xymatrix@R=52pt @C=55pt @!0{
&&& \VVV(\D_{3}, \C_{3}^{a})      &&&&& \VVV(\D_{3}, \C_{3}^{b}) &&\\
\VVV(\A_{4}) \ar@{-}[drr]& \VVV(\B_{4})\ar@{-}[dr] & & \VVV(\D_{3}, \D_{4})\ar@{-}[dl] \ar@{-}[dr]& \VVV(\C_{3}^{a}, \D_{4})  \ar@{-}[d] \ar@{-}[dr]& \VVV(\C_{4}^{a}) \ar@{-}[d]& \VVV(\C_{3}^{a}, \C_{3}^{b}) \ar@{-}[dl]\ar@{-}[drr]&   & \VVV(\C_{3}^{b}, \D_{4}) \ar@{-}[d]\ar@{-}[dllll] & \VVV(\C_{4}^{b}) \ar@{-}[dl]\\
&&\VVV(\D_{3}) \ar@/_-4pc/@{-}[uur]\ar@/_-7pc/@{-}[uurrrrrr] \ar@{-}[drr]&&\VVV(\D_{4}) \ar@{-}[d] & \VVV(\C_{3}^{a}) \ar@/_9pc/@{-}[uull]  \ar@{-}[dl]&&&  \VVV(\C_{3}^{b}) \ar@{-}[dllll] \ar@/_2.5pc/@{-}[uu]&\\
&&&& \VVV(\C_{2}) \ar@{-}[d]&&&&&\\
&&&& \textup{Trivial}  &&&&&
}
\]
\caption{Some varieties of positive S4-algebras}\label{Fig:Varieties}
\end{figure}
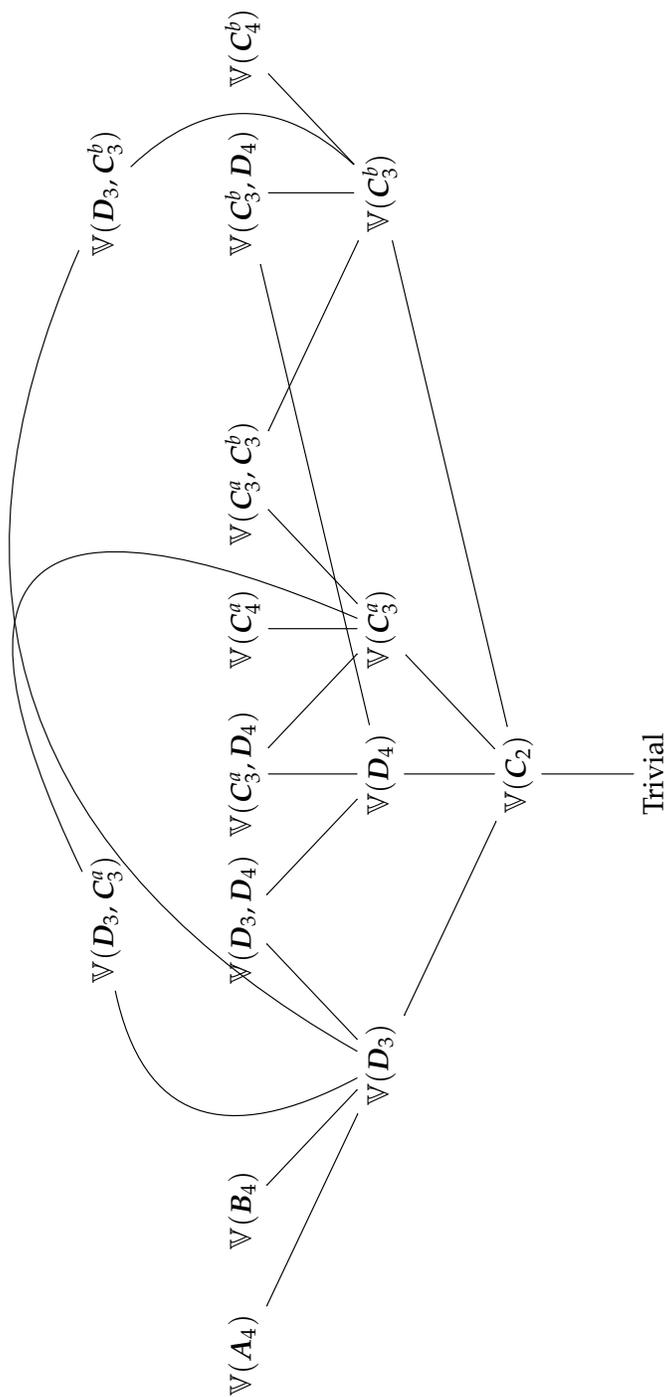
\end{landscape}

\bibliographystyle{plain}

\begin{thebibliography}{10}

\bibitem{BaKuVe15}
A.~Balan, A.~Kurz, and J.~Velebil.
\newblock Positive fragments of coalgebraic logics.
\newblock {\em {L}ogical {M}ethods in {C}omputer {S}cience}, 11(3):1--51, 2015.

\bibitem{Beklemishev14a}
L.~Beklemishev.
\newblock Positive provability logic for uniform reflection principles.
\newblock {\em {A}nnals of {P}ure and {A}pplied {L}ogic}, 165(1):82--105, 2014.

\bibitem{Bekle18}
L.~Beklemishev.
 A Note on Strictly Positive Logics and Word Rewriting Systems.
In {\em Larisa  Maksimova on Implication, Interpolation, and Definability},  volume~15 of {\em Outstanding Contributions to Logic}, pages 61--70.
\newblock {S}pringer-{V}erlag, 2018.

\bibitem{Be91b}
C.~Bergman.
\newblock Structural completeness in algebra and logic.
\newblock In {\em Algebraic logic ({B}udapest, 1988)}, volume~54 of {\em
  Colloq. Math. Soc. J\'anos Bolyai}, pages 59--73. North-Holland, Amsterdam,
  1991.

\bibitem{Be11g}
C.~Bergman.
\newblock {\em Universal Algebra: Fundamentals and Selected Topics}.
\newblock Chapman \& Hall Pure and Applied Mathematics. Chapman and Hall/CRC,
  2011.

\bibitem{BezGrigaml01}
G.~Bezhanishvili and R.~Grigolia.
\newblock {L}ocally tabular extensions of {MIPC}.
\newblock In M.~Zakharyaschev, K.~Segerberg, M.~de~Rijke, and H.~Wansing,
  editors, {\em {A}dvances in {M}odal {L}ogic}, volume~2, pages 101--120, 2001.

\bibitem{BlRiVe01}
P.~Blackburn, M.~de~Rijke, and Y.~Venema.
\newblock {\em Modal logic}.
\newblock Number~53 in Cambridge Tracts in Theoretical Computer Science.
  Cambridge University Press, Cambridge, 2001.


\bibitem{BP89}
W.~J. Blok and D.~Pigozzi.
\newblock {\em Algebraizable logics}, volume 396 of {\em Mem. Amer. Math. Soc.}
\newblock A.M.S., Providence, January 1989.

\bibitem{BuSa00}
S.~Burris and H.~P. Sankappanavar.
\newblock {\em A course in {U}niversal {A}lgebra}. The
  millennium edition, 2012.

\bibitem{CeJa97}
S.~A. Celani and R.~Jansana.
\newblock A new semantics for positive modal logic.
\newblock {\em Notre Dame Journal of Formal Logic}, 38(1):1--18, 1997.

\bibitem{CeJa99a}
S.~A. Celani and R.~Jansana.
\newblock Priestley duality, a {S}ahlqvist theorem and a {G}oldblatt-{T}homason
  theorem for positive modal logic.
\newblock {\em Logic Journal of the I.G.P.L.}, 7:683--715, 1999.

\bibitem{Ce06a}
S.~A. Celani.
\newblock Simple and subdirectly irreducibles bounded distributive lattices
  with unary operators.
\newblock {\em International Journal of Mathematics and Mathematical Sciences},
  2006:20, 2006.

\bibitem{ChZa97}
A.~Chagrov and M.~Zakharyaschev.
\newblock {\em Modal Logic}, volume~35 of {\em Oxford Logic Guides}.
\newblock Oxford University Press, 1997.

\bibitem{Citkin78a}
A.~Citkin.
\newblock On structurally complete superintuitionistic logics.
\newblock {\em Soviet Mathematics Doklady}, 19:816--819, 1978.

\bibitem{Cz01}
J.~Czelakowski.
\newblock {\em Protoalgebraic logics}, volume~10 of {\em Trends in
  Logic---Studia Logica Library}.
\newblock Kluwer Academic Publishers, Dordrecht, 2001.

\bibitem{CzDz}
J.~Czelakowski and W.~Dziobiak.
\newblock Congruence distributive quasivarieties whose finitely subdirectly
  irreducible members form a universal class.
\newblock {\em Algebra Universalis}, 27(1):128--149, 1990.

\bibitem{DaPr02}
B.~A. Davey and H.~A. Priestley.
\newblock {\em Introduction to lattices and order}.
\newblock Cambridge University Press, New York, second edition, 2002.

\bibitem{Day1975}
A.~Day.
\newblock Splitting algebras and a weak notion of projectivity.
\newblock {\em {A}lgebra {U}niversalis}, 5:153--162, 1975.

\bibitem{Du95}
J.~M. Dunn.
\newblock Positive modal logic.
\newblock {\em Studia Logica}, 55(2):301--317, 1995.

\bibitem{DzSt201x}
W.~Dzik and M.~Stronkowski.
\newblock Almost structural completeness; an algebraic approach.
\newblock {\em {A}nnals of {P}ure and {A}pplied {L}ogic}, 167(7):525--556, July
  2016.

\bibitem{DzWr73}
W.~Dzik and A.~Wro{\'n}ski.
\newblock Structural completeness of {G}\"odel's and {D}ummett's propositional
  calculi.
\newblock {\em Studia Logica}, 32:69--73, 1973.

\bibitem{AAL-AIT-f}
J.~M. Font.
\newblock {\em Abstract Algebraic Logic - An Introductory Textbook}, volume~60
  of {\em {S}tudies in {L}ogic - {M}athematical {L}ogic and {F}oundations}.
\newblock College Publications, London, 2016.

\bibitem{FJa09}
J.~M. Font and R.~Jansana.
\newblock {\em A general algebraic semantics for sentential logics}, volume~7
  of {\em Lecture Notes in Logic}.
\newblock A.S.L., second edition 2017 edition, 2009.

\bibitem{UAC11}
R.~Freese, E.~Kiss, and M.~Valeriote.
\newblock Universal {A}lgebra {C}alculator.
\newblock Available at: \url{www.uacalc.org}, 2011.

\bibitem{GeNaVe05}
M.~Gehrke, H.~Nagahashi, and Y.~Venema.
\newblock A {S}ahlqvist theorem for distributive modal logic.
\newblock {\em Annals of Pure and Applied Logic}, 131(1-3):65--102, 2005.

\bibitem{Gorbunov1976}
V.~A. Gorbunov.
\newblock Lattices of quasivarieties.
\newblock {\em {A}lgebra and {L}ogic}, 15:275--288, 1976.

\bibitem{Go98a}
V.~A. Gorbunov.
\newblock {\em Algebraic theory of quasivarieties}.
\newblock Siberian School of Algebra and Logic. Consultants Bureau, New York,
  1998.

\bibitem{HuCr96}
G.~E. Hughes and M.~J. Cresswell.
\newblock {\em A new introduction to modal logic}.
\newblock Routledge, London, 1996.

\bibitem{Ie01}
R.~Iemhoff.
\newblock On the admissible rules of intuitionistic propositional logic.
\newblock {\em The Journal of Symbolic Logic}, 66(1):281--294, 2001.

\bibitem{Ja02}
R.~Jansana.
\newblock Full models for positive modal logic.
\newblock {\em Mathematical Logic Quarterly}, 48(3):427--445, 2002.

\bibitem{Jerabek05c}
E.~Je{\v{r}}{\'a}bek.
\newblock Admissible rules of modal logic.
\newblock {\em {J}ournal of {L}ogic and {C}omputation}, 15(4):73--92, 2005.

\bibitem{Jo82}
P.~T. Johnstone.
\newblock {\em Stone Spaces}, volume~3 of {\em Cambridge studies in advanced
  mathematics}.
\newblock Cambridge University Press, 1982.

\bibitem{PTJo85}
P.~T. Johnstone.
\newblock Vietoris Locales and Localic Semilattices.
\newblock In {\em {C}ontinuous lattices and their {A}pplications}, volume 101, pages 155--180.
\newblock Marcel Dekker, 1985.

\bibitem{Jo70}
B.~J{\'o}nsson.
\newblock {\em Topics in universal algebra}, volume 250 of {\em Lecture Notes
  in Mathematics}.
\newblock Springer-Verlag, Berlin, 1970.

\bibitem{KoPi80}
P.~K{\"o}hler and D.~Pigozzi.
\newblock Varieties with equationally definable principal congruences.
\newblock {\em Algebra Universalis}, 11:213--219, 1980.

\bibitem{Kr99}
M.~Kracht.
\newblock {\em Tools and techniques in modal logic}, volume 142 of {\em Studies
  in Logic and the Foundations of Mathematics}.
\newblock North-Holland Publishing Co., Amsterdam, 1999.

\bibitem{MK07c}
M.~Kracht.
\newblock  Modal consequence relations.
\newblock In {\em {H}andbook of {M}odal {L}ogic}, volume~3, 491--545.
\newblock Elsevier Science Inc., New York, NY, USA, 2006.

\bibitem{MakRyb74}
L.~L. Maksimova and V.~V. Rybakov.
\newblock A lattice of normal modal logics.
\newblock {\em Algebra and Logic}, 13:105--122, 1974.

\bibitem{Ma71}
A.~I. Mal'cev.
\newblock {\em The metamathematics of algebraic systems, collected papers:
  1936-1967}.
\newblock Amsterdam, North-Holland Pub. Co., 1971.

\bibitem{McKe72}
R.~McKenzie.
\newblock Equational bases and nonmodular lattice varieties.
\newblock {\em Transactions of the Americal Mathematical Society}, 174:1--43,
  1972.

\bibitem{McMcTa87}
R.~N. McKenzie, G.~F. McNulty, and W.~F. Taylor.
\newblock {\em Algebras, lattices, varieties. {V}ol. {I}}.
\newblock The Wadsworth \& Brooks/Cole Mathematics Series. Wadsworth \&
  Brooks/Cole Advanced Books \& Software, Monterey, CA, 1987.

\bibitem{McKT44}
J.~C.~C. McKinsey and A.~Tarski.
\newblock The algebra of topology.
\newblock {\em Annals of Mathematics}, 45:141--191, 1944.

\bibitem{MeRo13a}
G.~Metcalfe and C.~R\"{o}thlisberger.
\newblock Admissibility in finitely generated quasivarieties.
\newblock {\em {L}ogical {M}ethods in {C}omputer {S}cience}, 9(2):1--19, 2013.

\bibitem{MorRafWan18}
T.~Moraschini, J.~G. Raftery, and J.~J. Wannenburg.
\newblock Singly generated quasivarieties and residuated structures.
\newblock Submitted manuscript, 2019.

\bibitem{Nguyen2000b}
L.~A. Nguyen.
\newblock {C}onstructing the {L}east {M}odels for {P}ositive {M}odal {L}ogic
  {P}rograms.
\newblock {\em {F}undamenta {I}nformaticae}, 42(2):29--60, 2000.

\bibitem{Nguyen05a}
L.~A. Nguyen.
\newblock On the complexity of fragments of modal logics.
\newblock In {K}ings~{C}ollege {P}ublications, editor, {\em {A}dvances in
  {M}odal {L}ogic}, volume~5, pages 249--268, 2005.

\bibitem{OlRaAl08}
J.~S. Olson, J.~G. Raftery, and C.~J. van Alten.
\newblock Structural completeness in substructural logics.
\newblock {\em Logic Journal of the I.G.P.L.}, 16(5):455--495, 2008.

\bibitem{Pa03}
A.~Palmigiano.
\newblock Coalgebraic semantics for positive modal logic.
\newblock In H.~P. Gumm, editor, {\em Electronic Notes in Theoretical Computer
  Science}, volume~82. Elsevier, 2003.

\bibitem{Pr70}
H.~A. Priestley.
\newblock Representation of distributive lattices by means of ordered {S}tone
  spaces.
\newblock {\em Bull. London Math. Soc.}, 2:186--190, 1970.

\bibitem{Pr72}
H.~A. Priestley.
\newblock Ordered topological spaces and the representation of distributive
  lattices.
\newblock {\em Proceedings of the London Mathematical Society. Third Series},
  24:507--530, 1972.

\bibitem{Ra06}
J.~G. Raftery.
\newblock Correspondences between {G}entzen and {H}ilbert systems.
\newblock {\em The Journal of Symbolic Logic}, 71(3):903--957, 2006.

\bibitem{RaxxNJ}
J.~G. Raftery.
\newblock Admissible {R}ules and the {L}eibniz {H}ierarchy.
\newblock {\em Notre {D}ame {J}ournal of {F}ormal {L}ogic}, 57(4):569--606,
  2016.

\bibitem{RaftSwiry16}
J.~G. Raftery and K.~{\'{S}}wirydowicz.
\newblock Structural completeness in relevance logics.
\newblock {\em Studia Logica}, 104(3):381--387, 2016.

\bibitem{Rauten1979bk}
W.~Rautenberg.
\newblock {\em Klassische und nichklassische Aussagenlogik}.
\newblock Vieweg Verlag, 1979.

\bibitem{ReV95-p}
J.~Rebagliato and V.~Verd{\'u}.
\newblock Algebraizable {G}entzen systems and the deduction theorem for
  {G}entzen systems.
\newblock Mathematics Preprint Series 175, University of Barcelona, 1995.


\bibitem{Blo77}
L. Rieger.
\newblock Zamerka o t.\ naz.\ svobodnyh algebrah s zamykanijami.
\newblock {\em Czechoslovak Mathematical Journal}, 7:16--20, 1957.

\bibitem{Ryb95}
V.~V. Rybakov.
\newblock Hereditarily structurally complete modals logics.
\newblock {\em The Journal of Symbolic Logic}, 60(1):266--288, March 1995.

\bibitem{Ry97}
V.~V. Rybakov.
\newblock {\em Admissibility of logical inference rules}, volume 136 of {\em
  Studies in Logic}.
\newblock Elsevier, Amsterdam, 1997.

\bibitem{Wro09}
A.~Wro{\'n}ski.
\newblock Overflow rules and a weakening of
  structural completeness.
  \newblock In {\em Rozwazania o {F}ilozfii {P}rawdziwej. {J}erezmu
  {P}erzanowskiemuw {Darze}}, 67--71.
\newblock Uniwersytetu Jagiello{\'n}skiego, Krak{\'o}w, 2009.

\end{thebibliography}

\end{document}